\documentclass[10pt,a4paper]{amsart}
\usepackage{palatcm}
\usepackage{amsthm,amssymb,amsmath}

\allowdisplaybreaks

\def\phi{\varphi}
\newcommand{\D}{\mathcal{D}}
\newcommand{\B}{\mathcal{B}}
\newcommand{\M}{\mathcal{M}}
\newcommand{\T}{\mathcal{T}}
\newcommand{\R}{\mathbb{R}}
\newcommand{\C}{\mathbb{C}}
\newcommand{\CC}{\mathcal{C}}
\newcommand{\CCC}{\mathfrak{C}}
\newcommand{\Dg}{\mathcal{D}\gamma}
\newcommand{\dW}{d\mathrm{W}(\gamma)}
\def\P{\mathfrak{P}}
\newcommand{\Pcl}{\mathfrak{P_{cl}}}
\def\oooo{\mathrm{o}}
\def\OOOO{\mathrm{O}}
\newcommand{\Hom}{\mathrm{Hom}}
\newcommand{\vol}{\mathrm{vol}}
\newcommand{\hol}{\mathrm{hol}}
\newcommand{\scal}{\mathrm{scal}}
\newcommand{\ric}{\mathrm{ric}}
\newcommand{\Ric}{\mathrm{Ric}}
\newcommand{\tr}{\mathrm{tr}}
\newcommand{\Tr}{\mathrm{Tr}}
\newcommand{\grad}{\mathrm{grad}}
\def\div{\mathrm{div}}
\def\Laenge{\mathrm{L}}
\newcommand{\LL}{\mathcal{L}}
\newcommand{\E}{\mathrm{E}}
\newcommand{\injrad}{\mathrm{injrad}}
\newcommand{\id}{\mathrm{id}}
\def\epsilon{\varepsilon}
\newcommand{\eps}{\varepsilon}
\newcommand{\<}{\langle}
\def\>{\rangle}

\newcommand{\gabxy}{\gamma_{x,y}^{a,b}}
\newcommand{\gotxy}{\gamma_{x,y}^{0,t}}
\newcommand{\hvt}{\widehat v_t}
\newcommand{\hVt}{\widehat V_t}
\newcommand{\tvt}{\widetilde v_t}
\newcommand{\tVt}{\widetilde V_t}
\newcommand{\bvt}{\bar v_t}
\newcommand{\hwt}{\widehat w_t}
\newcommand{\hWt}{\widehat W_t}
\newcommand{\hW}{\widehat W}
\newcommand{\hw}{\widehat w}

\newcommand{\abl}[1]{\frac{\partial}{\partial #1}}

\theoremstyle{plain}
\newtheorem{thm}{Theorem}[section]

\newtheorem{cor}[thm]{Corollary}
\newtheorem{lem}[thm]{Lemma}

\theoremstyle{definition}
\newtheorem{dfn}[thm]{Definition}
\newtheorem{rem}[thm]{Remark}
\newtheorem{rems}[thm]{Remarks}
\newtheorem{exm}[thm]{Example}

\begin{document}

\title{{Path integrals on manifolds by finite dimensional approximation}}

\author{Christian B\"ar and Frank Pf\"aff\-le}

\address{Universit\"at Potsdam\\
Institut f\"ur Mathematik\\
Am Neuen Palais 10\\
14469 Potsdam\\
Germany}

\email{baer@math.uni-potsdam.de, pfaeffle@math.uni-potsdam.de}

\subjclass[2000]{58J65, 58J35, 47D06}
\keywords{path integral, Feynman-Kac formula, Chernoff's theorem, generalized
  Laplace operator, Riemannian manifold, heat equation, heat kernel,
  Hess-Schrader-Uhlenbrock estimate}

\date{\today}

\maketitle

\begin{abstract}
Let $M$ be a compact Riemannian manifold without boundary and let $H$ be a
self-adjoint generalized Laplace operator acting on sections in a bundle over
$M$.
We give a path integral formula for the solution to the corresponding heat
equation.
This is based on approximating path space by finite dimensional spaces of
geodesic polygons.
We also show a uniform convergence result for the heat kernels.
This yields a simple and natural proof for the Hess-Schrader-Uhlenbrock
estimate and a path integral formula for the trace of the heat operator.
\end{abstract}

\setcounter{tocdepth}{1}
\tableofcontents

\section{Introduction}

Many diffusion processes such as heat flow or gas diffusion are mathematically
desribed by the {\em heat equation}
\begin{equation}
\frac{\partial U}{\partial t} + HU =0.
\label{eq:Heat}
\end{equation}
The operator $H$ is typically of the form $H=\Delta +V$ where $\Delta$ ist the
Laplace operator and $V$ some potential.
We use the sign convention that $H$ be spectrally bounded from below, i.~e.\
for the Laplace operator on $\R^m$ we have $\Delta = - \sum_{j=1}^m
\frac{\partial^2}{\partial (x^j)^2}$.
Motivated by the microscopic picture of diffusion one expects the solution $U$
of (\ref{eq:Heat}) with given initial condition $U(0,x)=u(x)$ to be given by a
``path integral'' of the form
\begin{equation}
U(t,x) = 
\frac{1}{Z} \, \int_{\CCC_x(\R^m,t)} 
\exp\left(-\frac12\E(\gamma) - \int_{0}^{t} V(\gamma(s))\, ds\right)
\cdot u(\gamma(t))\,\Dg .
\label{eq:Path}
\end{equation}
Here $\CCC_x(\R^m,t)$ is the set of all continuous paths $\gamma$ in $\R^m$
parametrized on $[0,t]$ starting at $\gamma(0)=x$.
Moreover, $\E(\gamma)$ denotes the energy of $\gamma$, $\Dg$ is a suitable
measure on $\CCC_x(\R^m,t)$, and $Z$ is a normalizing constant.

Problems with a too naive approach to formula (\ref{eq:Path}) are numerous:
the space $\CCC_x(\R^m,t)$ is infinite dimensional and the measure $\Dg$ does not
exist, the energy is defined only for differentiable paths which are expected
not to contribute to the integral, and the normalizing constant $Z$ is
infinite.
It turns out that the measure $\dW := \frac1Z
\exp\left(-\frac12\E(\gamma)\right) \Dg$ does exist and under reasonable
conditions 
\begin{equation}
U(t,x) = 
 \int_{\CCC_x(\R^m,t)} 
\exp\left(- \int_{0}^{t} V(\gamma(s))\, ds\right)
\cdot u(\gamma(t))\,\dW 
\label{eq:FeynmanKac}
\end{equation}
holds true.
This is known as the {\em Feynman-Kac formula} and $\dW$ is called {\em Wiener
measure}.
There is a huge literature on these topics, see e.~g.\ \cite{S} and \cite{JL}
and the references therein.

After replacing $\R^m$ by a closed Riemannian manifold $M$ formula
(\ref{eq:FeynmanKac}) still holds.
One should note however that often Wiener measure on a manifold is defined in
such a way that (\ref{eq:FeynmanKac}) becomes tautological in case $V=0$.
The Feynman-Kac formula for non-trivial $V$ is then a rather simple
consequence of the Trotter formula.
In \cite{AD} the authors approximate path space $\CCC_x(M,T)$ by finite
dimensional spaces of geodesic polygons and obtain two approximations for
Wiener measure.
They differ by a scalar curvature term.
We will come back to this.
Corollary 1.9 in \cite{AD} also gives a non-tautological path integral for
$H=\Delta$ by finite dimensional approximation.

In the present paper we are concerned with general self-adjoint Laplace type
operators acting on sections in a vector bundle over $M$.
So we allow for systems of equations rather than scalar equations only.
Our main result is Theorem~\ref{thm:PathIntegral1} where we give a path
integral expression for solutions to the heat equation for such general
operators by finite dimensional approximation.
In case the potential $V$ is scalar valued one could state the result formally
as 
\begin{eqnarray}
\lefteqn{U(t,x)\quad =}\nonumber\\ 
&& 
\frac{1}{Z} \, \int_{\CCC_x(M,t)} 
\exp\left(-\frac12\E(\gamma) 
+ \int_{0}^{t} \left(\frac13\scal(\gamma(s))-V(\gamma(s))\right) ds\right)
\cdot \tau(\gamma)_t^0\cdot u(\gamma(t))\,\Dg ,
\label{eq:PathGeneral}
\end{eqnarray}
compare Corollary~\ref{cor:PathIntegral1}.
Here $\tau(\gamma)$ denotes parallel translation along $\gamma$.
There are stochastic versions of this path integral using Wiener measure and
stochastic parallel transport, see e.~g.\ \cite[Sec.~4.1]{DT}.

Our technique allows us to derive different versions of the path integral
formula.
For example, one can remove the scalar curvature term in
(\ref{eq:PathGeneral}) if one uses another measure on the approximating spaces
of geodesic polygons.
This was remarked already in \cite{AD}.
We show that one can actually interpolate between the two formulas, compare
Theorem~\ref{thm:PathIntegral3}.

Theorem~\ref{thm:PathKernel1} says that for suitable approximations even the
heat kernels, i.~e., the integral kernels of the solution operators for the
heat equation, are uniformly approximated by the corresponding kernels
obtained by integration over the spaces of geodesic polygons.
As a consequence we find a very simple and natural proof of the
Hess-Schrader-Uhlenbrock estimate for the heat kernel by the kernel of a
scalar comparison operator.
Moreover, we can express the trace of the heat operators by a path integral.
Formally, Theorem~\ref{thm:Trace} says in case $V$ is scalar
\begin{eqnarray*}
\lefteqn{\Tr(e^{-tH}) \quad =}\\
&&
\frac{1}{Z} \,
\int_{\CCC_{\mathfrak{cl}}(M,t)} 
\exp\left(-\frac12\E(\gamma)
+\int_0^{t}\left(\frac13 \scal(\gamma(s))-V(\gamma(s))\right)\, 
ds\right)\tr(\hol(\gamma))\,\Dg .
\end{eqnarray*}
Here $\CCC_{\mathfrak{cl}}(M,t)$ denotes the space of closed continuous loops
in $M$, parametrized on $[0,t]$, and $\hol(\gamma)$ is the holonomy of such a
loop $\gamma$.

Despite the heuristic formal expressions we have given in this introduction
all results of the paper involve well-defined quantities only
and the derivations are mathematically rigorous.
Our technique of proof is based on the concept of Chernoff equivalence of
families of operators in a version we learnt from \cite{SWW07}.
This allows one to organize the analysis in a quite transparent manner.
The short time asymptotics of the heat kernel also play an important role.

The present article might be regarded as a sequel to \cite{AD} but we do not
use any of the results therein.
The derivation of our general path integral formulas is entirely
self-contained.
In \cite[Thm.~3]{SWW07} quite general results on Feller semigroups
are applied to conclude from path integral formulas like (\ref{eq:PathGeneral})
that the associated measures on path space converge weakly to the Wiener
measure. 
This reproves the approximation of Wiener measure in \cite{AD}.
The same methods would apply to our results in Section~\ref{sec:Alternative}
to obtain other versions of this result.

In the stochastically oriented literature it is common to look at $H = \frac12
\Delta +V$ while we use the convention $H = \Delta +V$.
This is the reason why we have a factor $\frac13$ in (\ref{eq:PathGeneral})
where other works like \cite{AD} have $\frac16$.

{\bf Acknowledgements.}
It is a pleasure to thank M.~Klein and S.~Roelly for helpful discussion
and SFB 647 for financial support.

\section{Chernoff's theorem}\label{sec:cherno}

In this section we give some definitions and collect some approximation
results for semigroups of bounded 
operators, in particular Chernoff's Theorem in the version of
\cite[Section~2]{SWW07}.
There is no claim of originality for this section.
We include it for the convenience of the reader and to set up notation.

Throughout this section let $X$ and $Y$ be Banach spaces. 
Denote $\LL(X,Y)$ the space of bounded linear operators from $X$ to $Y$, and 
let $\B(X)=\LL(X,X)$ be the space of bounded linear operators on $X$.
 
\begin{dfn}
A family of operators $T_i\in\LL(X,Y)$, $i\in I$, is called {\em equicontinuous}
if the operator norms are uniformly bounded,
\[ \sup_{i\in I} \|T_i \|<\infty. \]
\end{dfn}

\begin{thm}[Banach-Steinhaus]
Let $T_i\in\LL(X,Y)$, $i\in I$, such that
\[\sup_{i\in I}\|T_iu\| <\infty\]
for each $u\in X$. 
Then the family $(T_i)_{i\in I}$ is equicontinuous.
\end{thm}
For a proof see e.~g.\ \cite[Thm.~III.9]{RS}.

\begin{lem}\label{lem:compactuniform}
Let $(T_n)_{n\ge 1}\subset \LL(X,Y)$ be an equicontinuous sequence of
operators, and let $K\subset X$ be a compact subset such that
\[ T_nu\xrightarrow{n\to\infty}0 \;\mbox{ for each }u\in K.\]
Then one has uniform convergence on $K$,
\[ \lim_{n\to\infty} \sup_{u\in K} \|T_nu \|=0. \]
\end{lem}
\begin{proof}
Suppose the opposite is true.
Then there is an $\varepsilon >0$ such that for any $n\ge 1$ there exists a
  $u_n\in K$ with $\|T_nu_n \|\ge\varepsilon$.
Since $K$ is compact $(u_n)_{n\ge1}$ converges, after passing to a subsequence,
to some $u_\infty\in K$.
By equicontinuity one gets a constant $C>0$ with
\[ \|T_n(u_\infty-u_n)\|\le C\cdot\|u_\infty-u_n
\|\xrightarrow{n\to\infty} 0.\]
The inverse triangle inequality yields
\[
\|T_nu_\infty\|=\|T_n(u_\infty-u_n)+T_nu_n\|\ge
\Big|\|T_n(u_\infty-u_n)\|-\|T_nu_n\| \Big| \ge\frac{\varepsilon}{2}
\]
for large $n$, in contradiction to $T_nu_\infty\xrightarrow{n\to\infty} 0$.  
\end{proof}

\begin{lem}\label{lem:equidense}
Let $(T_n)_{n\ge 1}$ be a sequence of operators in $\LL(X,Y)$, let
$T_\infty\in \LL(X,Y)$, and let $\D\subset X$ be a dense subset. 
Then the following statements are equivalent:
\begin{enumerate}
\item[a)]
$(T_n)_{n\ge 1}$ converges strongly to $T_\infty$, i.~e., 
$\lim_{n\to\infty} T_nu=T_\infty u$ for all $u\in X$.
\item[b)]
$(T_n)_{n\ge 1}$ is equicontinuous and $\lim_{n\to\infty} T_nu=T_\infty u$ for
  all $u\in \D$. 
\end{enumerate}
\end{lem}

The conclusion ``a) $\Rightarrow$ b)'' is a direct consequence of the
Banach-Steinhaus Theorem.
The opposite conclusion is elementary.

\begin{dfn}\label{proper}
A map $S:[0,\infty)\to \B(X)$ is called a {\em proper family} if the following
  holds
\begin{enumerate}
\item[a)] $S_0=\id_X$,
\item[b)] $S$ is strongly continuous, i.~e., for any $u\in X$ the map
  $[0,\infty)\to X$ given by $t\mapsto S_tu$ 
  is continuous,
\item[c)] $\|S_t\|=1 + \OOOO(t)$ as $t\searrow 0$, and 
\item[d)] there exists a (possibly unbounded) closed operator $L$ in $X$ with
dense domain $\D(L)$ which is the generator of a strongly continuous semigroup  
$\left(e^{tL}\right)_{t\ge 0}$ on $X$ such that
\[
\frac{S_t-\id_X}{t}u\xrightarrow{t \searrow 0} Lu
\]
for all $u\in X$ of the form $u=e^{aL}v$ with $a>0$ and $v\in \D(L)$.
This operator $L$ will be denoted by $DS$.
\end{enumerate}
\end{dfn}

\begin{rems}\label{rmk1}
a) 
Since the operators $e^{aL}$, $a>0$, leave the domain $\D(L)$ invariant we
have $\D'(L):=\{e^{aL}v\,|\,v\in\D(L),\,a>0\} \subset \D(L)$.
Moreover, since $\D(L)$ is dense in $X$ and since $e^{aL}v \to v$ as
$a\searrow 0$ the subset $\D'(L)$ is dense in $X$ as well.\\
b) 
Any strongly continuous semigroup forms a proper family.\\
c) 
Theorem~\ref{thm:Chernoff} below shows that one can reconstruct the
semigroup $\left(e^{tL}\right)_{t\ge 0}$ and hence also its generator $L$
from the proper family $S$. 
Therefore the operator $L$ is uniquely determined by $S$ and the notation
$L=DS$ is justified.
\end{rems}

\begin{dfn}
A tuple of positive real numbers $\T=(t_1,\ldots,t_r)$ is called a
\emph{partition}.
Its \emph{length} is given by $\Laenge(\T):=t_1+\ldots+t_r$ and its \emph{mesh} by
$|\T|:=\max_{j=1,\ldots,r}t_j$.
\end{dfn}

We think of $\T$ as a subdivision of the interval $[0,t]$, $t=\Laenge(\T)$, into $r$
subintervals $[0,t_1]$, $[t_1,t_1+t_2]$, $\ldots\,$, $[t-t_r,t]$. 

Now we show a version of Chernoff's Theorem following \cite[Prop.~1]{SWW07}
and \cite[Prop.~3]{SWW03}. 
The main technical advantage of this version as compared to the original
version of Chernoff's Theorem (see \cite{Ch}) consists of the fact that we
need not assume the $S_t$ to be contractions and the partitions $\T_n$ need not
be equidistant.

\begin{thm}\label{thm:Chernoff}
Let $(S_t)_{t\geq 0}$ be a proper family of bounded linear operators on $X$
with $L=DS$. 
Consider partitions $\T_n=(t_1^n,\ldots,t_{r_n}^n)$, $n\ge 1$, with
$\Laenge(\T_n)\to\tau > 0$ and $|\T_n|\to 0$ as $n\to\infty$.
Then for any $u\in X$ one has
\begin{equation}\label{eq:Chernoff}
 S_{t_1^n}\cdots S_{t_{r_n}^n}u \xrightarrow{n\to\infty} e^{\tau L}u. 
\end{equation}
\end{thm}
\begin{proof}
By c) in Definition~\ref{proper} there exists $q>0$ such that $\|S_t\|\le
e^{qt}$ for all sufficiently small $t\ge 0$.
Now $|\T_n|\to 0$ as $n\to\infty$, hence for sufficiently large $n$ one gets 
\begin{equation}\label{eq:Stsexp}
 \left\|S_{t_1^n}\cdots S_{t_{j}^n}\right\| 
 \leq e^{qt_1^n}\cdots e^{qt_j^n}
 \le e^{q\cdot \Laenge(\T_n)}
\end{equation}
for all $j=1,\ldots,r_n$.

Fix $a>0$.
We note that $\D(L)$ equipped with the graph norm $\|u\|_L = \|u\| + \|Lu\|$
is a Banach space.
With respect to this norm the operator  
\[
\left( \frac{1}{t}\left(S_t-\id_X\right)-L \right)e^{aL}:\D(L)\to X
\]
is bounded for any $t>0$.
By b) and d) in Definition~\ref{proper} we get
$ \sup_{0< t\le t'}
\left\| \left( \frac{1}{t}\left(S_t-\id_X\right)-L \right)e^{aL}v
\right\| <\infty $
for any $v\in \D(L)$ and any $t'>0$.
Hence, by the Banach-Steinhaus theorem, the family 
$\left\{ \left(
\frac{1}{t}\left(S_t-\id_X\right)-L \right)e^{aL}\right\}_{0< t\le t'}$ is
equicontinuous on $\D(L)$. 

For each $v\in\D(L)$ the map $[0,\infty)\to\D(L)$, $s\mapsto e^{sL}v$, is
  continuous (w.~r.~t.\ the graph norm on $\D(L)$).
Therefore $\left\{e^{sL}v \,|\,s\in[0, b-a]\right\}$ is a compact subset of
  $\D(L)$ for any $v\in\D(L)$ and any $b>a$.
Lemma~\ref{lem:compactuniform} and d) in Definition~\ref{proper} then give
\begin{equation}\label{eq:glmlim}
\lim_{t\to 0} \sup_{a\le s\le b} \left\| 
\left( \frac{1}{t}\left(S_t-\id_X\right)-L \right)e^{sL}v \right\|
=
\lim_{t\to 0} \sup_{0\le s\le b-a} \left\| 
\left( \frac{1}{t}\left(S_t-\id_X\right)-L \right)e^{aL}e^{sL}v \right\|
=0.
\end{equation}
For each partition $\T_n=(t_1^n,\ldots,t_{r_n}^n)$ we put
$s_j^n:=\sum_{i=j+1}^{r_n}  t_i^n$ where
$j=1,\ldots,r_n$,
and we get
\[
S_{t_1^n}\cdots S_{t_{r_n}^n}-e^{s_0^n\cdot L}=\sum_{j=1}^{r_n} S_{t_1^n}\cdots
S_{t_{j-1}^n} \cdot\left(S_{t_j^n}-e^{t_j^n\cdot L} \right)\cdot e^{s_j^n\cdot
  L}. 
\]
By (\ref{eq:Stsexp}) and by $\Laenge(\T_n)\to\tau$ there is a number $C>0$ 
such that for all $n$ one has
\begin{equation}\label{eq:Stglmst}
\left\| S_{t_1^n}\cdots S_{t_j^n}\right\|\le C
\quad\mbox{ for all } j=1,\ldots,r_n.
\end{equation}
For any $u\in X$ this yields 
\begin{eqnarray}
\lefteqn{\left\|\left(S_{t_1^n}\cdots S_{t_{r_n}^n}-e^{s_0^n\cdot L}\right)
    u\right\|} \nonumber\\
&\le& 
C\cdot\sum_{j=1}^{r_n}\left\| \left(S_{t_j^n}-e^{t_j^n\cdot L} \right)
\cdot e^{s_j^n \cdot L}u \right\| \nonumber\\
&\le& 
C\cdot\sum_{j=1}^{r_n} t_j^n\cdot \left\| 
\left(\frac{S_{t_j^n}-\id_X}{t_j^n}-\frac{e^{t_j^n\cdot
    L}-\id_X}{t_j^n}\right) 
\cdot e^{s_j^n\cdot L}u \right\|\nonumber\\
&\le& 
C\cdot L(\T_n)\cdot \sup_{\substack{0<t\le |\T_n| \\ 0\leq s\le L(\T_n)}}
\left\| \left(\frac{S_{t}-\id_X}{t}-\frac{e^{t\cdot L}-\id_X}{t}\right)
\cdot e^{s\cdot L}u \right\| .
\label{est:supmu}
\end{eqnarray}
To handle the right hand side of (\ref{est:supmu}) we apply (\ref{eq:glmlim})
twice, once for the proper family $(S_t)_{t\ge 0}$ and once for the proper
family $(e^{t L})_{t\ge 0}$.
For $u$ of the form $u=e^{aL}v$ with $v\in\D(L)$ and $a>0$ we get
\begin{eqnarray*}
0
&\le& 
\lim_{n\to\infty}\sup_{\substack{0< t\le |\T_n| \\ 0\leq s\le L(\T_n)}}
\left\| \left(\frac{S_{t}-\id_X}{t}-\frac{e^{t\cdot L}-\id_X}{t}\right)
\cdot e^{s\cdot L}u \right\|\\
&=& 
\lim_{n\to\infty}\sup_{\substack{0< t\le |\T_n| \\ 0\leq s\le L(\T_n)}}
\left\| \left(\frac{S_{t}-\id_X}{t}-\frac{e^{t\cdot L}-\id_X}{t}\right)
\cdot e^{(s+a)\cdot L}v \right\|\\
&\le& 
\lim_{n\to\infty}\sup_{\substack{0< t\le |\T_n| \\ a\leq s\le b}}
\left\| \left(\frac{S_{t}-\id_X}{t}-\frac{e^{t\cdot L}-\id_X}{t}\right)
\cdot e^{s\cdot L}v \right\|\\
&\le& 
\lim_{t\to 0}
\left( 
\sup_{a\le s\le b}\left\| 
\left( \frac{1}{t}\left(S_t-\id_X\right)-L \right)e^{sL}v \right\| +
\sup_{a\le s\le b}\left\| 
\left( \frac{1}{t}\left(e^{t L}-\id_X\right)-L \right)e^{sL}v \right\|
\right)\\
&\stackrel{(\ref{eq:glmlim})}{=}& 0
\end{eqnarray*}
where $b>0$ is chosen such that $L(\T_n)-a\le b$ for all $n$.  
Thus
$$
\|(S_{t_1^n}\cdots S_{t_{r_n}^n} - e^{s_0^n L})u\| 
\xrightarrow{n\to\infty} 0 .
$$
Now we notice that $s_0^n=L(\T_n)\to\tau$ for $n\to\infty$
and hence $e^{s_0^n\cdot L}u\to e^{\tau L}u$ for any $u\in X$.
Therefore we have shown the claim (\ref{eq:Chernoff}) for $u=e^{aL}v$ with
$a>0$ and $v\in\D(L)$,
\[
S_{t_1^n}\cdots S_{t_{r_n}^n}e^{a L}v \xrightarrow{n\to\infty} e^{\tau
 L} e^{aL}v.
\]
We observe that $\left\{e^{aL}v\,\mid\, a>0,v\in\D(L) \right\}\subset X$ is
dense by Remark~\ref{rmk1} a) and that the operators $T_n:=S_{t_1^n}\cdots
S_{t_{r_n}^n}$ form an equicontinuous sequence $(T_n)_{n\ge 1}$ by 
(\ref{eq:Stglmst}).
Applying Lemma~\ref{lem:equidense} concludes the proof. 
\end{proof}

\begin{dfn}
Two proper families $S=(S_t)_{t\ge 0}$ and $T=(T_t)_{t\ge 0}$ are
called \emph{Chernoff equivalent} if $DS=DT$.
\end{dfn}

Two Chernoff equivalent families yield the same semigroup via the
approximation (\ref{eq:Chernoff}) in Theorem~\ref{thm:Chernoff}.
The next two lemmas give us sufficient criteria for Chernoff equivalence.

\begin{lem}\label{lem:ChernoffKriterium}
Let $X$ be a Banach space, let $[0,\infty) \to \B(X)$, $t \mapsto S_t$, be a
proper family.
Moreover, let $(0,\infty) \to \B(X)$, $t \mapsto T_t$, be a strongly continuous
family such that 
$$
\|S_t-T_t\| = \oooo(t) \mbox{ as } t \searrow 0.
$$
Then putting $T_0:=\id_X$ yields a proper family $[0,\infty) \to \B(X)$, $t
\mapsto T_t$, which is Chernoff equivalent to the family $(S_t)_{t\geq0}$.
\end{lem}

\begin{proof}
For fixed $u\in X$ we have $\|T_tu-u\| \leq \|T_t-S_t\|\cdot\|u\| + \|S_tu-u\|
\to 0$ as $t\searrow 0$.
Hence putting $T_0:=\id_X$ yields an extension of the family which is strongly
continuous on $[0,\infty)$.

From $|\|T_t\|-\|S_t\|| \leq \|T_t-S_t\| =\oooo(t)$ and $\|S_t\|=1+\OOOO(t)$ we
conclude $\|T_t\|=1+\OOOO(t)$.

Finally, for any $u\in\D'(L)$, $L=DS$, we have 
$$
\left\| \frac{T_t-\id_X}{t}u-Lu\right\| 
\leq 
\left\| \frac{S_t-\id_X}{t}u-Lu\right\| + \frac1t \|T_t-S_t\|\cdot\|u\| \to 0 
$$
as $t\searrow 0$.
This proves the claim.
\end{proof}

The following variation of the criterion for Chernoff equivalence can be found
in \cite[Lemma~1]{SWW07}.

\begin{lem}\label{lem:ChernoffKriterium1}
Let $X$ be a Banach space, let $(S_t)_{t\geq 0}$, be a proper family of
bounded operators on $X$ with $L=DS$.
Moreover, let $(T_t)_{t>0}$, be a strongly continuous family of bounded
operators on $X$ such that 
$$
\|T_t\| = 1+ \OOOO(t) \mbox{ as } t \searrow 0
$$
and 
$$
\|S_tu-T_tu\| = \oooo(t) \mbox{ as } t \searrow 0
$$
for all $u\in X$ of the form $u=e^{aL}v$ where $a>0$ and $v\in X$.

Then putting $T_0:=\id_X$ yields a proper family $(T_t)_{t\geq 0}$ which is
Chernoff equivalent to the family $(S_t)_{t\geq0}$. 
\end{lem}

\begin{proof}
For $u$ of the form $u=e^{aL}v$ with $a>0$ we have
$$
\|T_tu - u \|
\quad\leq\quad
\|T_tu-S_tu\| + \|S_tu-u\|
\quad=\quad
\oooo(t) + \oooo(1)
\quad=\quad
\oooo(1).
$$
Thus $\|T_tu - u \| \xrightarrow{t\searrow 0} 0$ for all $u$ in a dense subset
of $X$.
From $\|T_t\|=1+\OOOO(t)$ we see that the family $(T_t)_{0\leq t\leq
  \varepsilon}$ is equicontinuous.
By Lemma~\ref{lem:equidense} we have that $t\mapsto T_t$ is strongly
continuous also at $t=0$.
Hence conditions a), b), and c) in Definition~\ref{proper} are satisfied for
the family $(T_t)_{t\geq 0}$.

For $u=e^{aL}v$ with $v\in \D(L)$ we have
$$
\left\| \frac{T_t-\id_X}{t}u-Lu\right\| 
\quad\leq\quad 
\left\| \frac{S_t-\id_X}{t}u-Lu\right\| + \frac1t
\underbrace{\|T_tu-S_tu\|}_{=\oooo(t)} 
\quad\to\quad 0  
$$
as $t\searrow 0$.
This proves the lemma.
\end{proof}

\section{The main theorem}

\textbf{General assumption.}
From now on we adopt the following notation.
Let $M$ be an $m$-dimensional compact Riemannian manifold without boundary,
let $E\to M$ be a real or complex vector bundle, equipped with a Riemannian or
Hermitian metric $\<\cdot,\cdot\>$ respectively.
For the necessary basics on Riemannian geometry we refer the reader to
\cite{C}. 

By $H$ we denote a formally self-adjoint generalized Laplace operator 
acting on sections in $E$.
Here ``generalized Laplace operator'' means that the principal symbol of $H$
is given by the Riemannian metric of $M$, i.~e.\ in local coordinates we have, 
$$
H = -\sum_{j,k=1}^m g^{jk}\frac{\partial^2}{\partial x^j \partial x^k} 
+ \mbox{ lower order terms}.
$$
We assume that $H$ has smooth coefficients.
Formal self-adjointness means that for all smooth sections $u$ and $v$ in $E$
$$
(Hu,v) = (u,Hv)
$$
holds where $(u,v) = \int_M \< u(x),v(x)\>\,dx$ is the corresponding
$L^2$-scalar procuct.
Here $dx$ denotes the volume measure induced by the Riemannian metric.
It is well-known that $H$ is essentially self-adjoint in the Hilbert space
$L^2(M,E)$ of square-integrable sections in $E$ when given the domain
$C^\infty(M,E)$ of smooth sections in $E$, see e.~g.\ \cite[Prop.~2.33,
p.~89]{BGV}. 
Moreover, one knows that $H$ can be written in the form 
$$
H = \nabla^*\nabla + V
$$
where $\nabla$ is a metric connection on $E$ and $V$ is a smooth section in
symmetric endomorphisms of $E$, compare \cite[Prop.~2.5, p.~67]{BGV}.
We call $\nabla$ the \emph{connection determined by $H$} and $V$ its
\emph{potential}. 

\begin{exm}
The simplest example for $H$ as described above is the \emph{Laplace-Beltrami
operator} $H=\Delta$ acting on functions.
Here $E$ is the trivial real line bundle, $\nabla = d$ the usual derivative
and $V=0$.
\end{exm}

\begin{exm}
More generally, let $E=\bigwedge^k T^*M$ be the bundle of $k$-forms.
Then we may take the \emph{Hodge Laplacian} $H=d\delta + \delta d$ acting on
$k$-forms.
Here $d$ denotes exterior differentiation and $\delta$ its formal adjoint.
The Weitzenb\"ock formula says that $H = \nabla^*\nabla + V$ where
$\nabla$ is the Levi-Civita connection and $V$ depends linearly on the
curvature tensor of $M$. 
For example, for $k=1$ we have $V=\Ric$, see e.~g.\ \cite[Ch.~1.I]{Be}.
\end{exm}

\begin{exm}
If $M$ is a spin manifold one can form the spinor bundle $E=\Sigma M$ and the
Dirac operator $D$ acting on sections in $E$.
Then $H=D^2 = \nabla^*\nabla + \frac14 \scal$ is a self-adjoint generalized
Laplace operator.

More generally, the square of any generalized Dirac operator in the sense of
Gromov and Lawson yields a self-adjoint generalized Laplacian, see e.~g.\
\cite[Sec.~1,2]{GL}.
\end{exm}

By functional calculus the self-adjoint extension of $H$ generates a
strongly continuous semigroup $t\mapsto e^{-tH}$ in the Hilbert space
$L^2(M,E)$. 
For $u \in L^2(M,E)$ the section $U(t,x) := (e^{-tH}u)(x)$, $(t,x) \in
[0,\infty)\times M$, is the unique solution to the heat equation
$$
\frac{\partial U}{\partial t} + HU =0
$$
satisfying the inition condition $U(0,x)=u(x)$.
The aim of this article is to derive a ``path integral formula'' for $U$.

\begin{dfn}
Let $\T=(t_1,\ldots,t_r)$ be a partition of length $t$.
Put $\sigma_j(\T) := t_1 + \ldots + t_j$.
Let $\gamma : [0,t] \to M$ be a continuous curve.
Put $x_j := \gamma(\sigma_j(\T))$.

The curve $\gamma$ is called a {\em geodesic polygon} in $M$
with respect to $\T$ if any two subsequent points $x_j$ and $x_{j+1}$ are not 
cut-points of each other and $\gamma|_{[\sigma_j(\T),\sigma_{j+1}(\T)]}$ is
the unique shortest geodesic joining them.
\end{dfn}

The assumption on $x_j$ and $x_{j+1}$ of not being cut-points is made to
ensure uniqueness of the shortest geodesic joining $x_j$ and $x_{j+1}$.
Note in particular, that when restricted to one of the subintervals
$[\sigma_j(\T),\sigma_{j+1}(\T)]$, the curve $\gamma$ is smooth and
parametrized proportionally to arc-length.

We denote the set of all geodesic polygons in $M$ with respect to $\T$ by
$\P(M,\T)$.
For $x,y\in M$ we set $\P_x(M,\T) := \{\gamma\in\P(M,\T)\,|\, \gamma(0)=x\}$,
$\P^y(M,\T) := \{\gamma\in\P(M,\T)\,|\, \gamma(t)=y\}$, and $\P_x^y(M,\T) :=
\{\gamma\in\P(M,\T)\,|\, \gamma(0)=x \mbox{ and }\gamma(t)=y\} = \P_x(M,\T)
\cap \P^y(M,\T)$.
Moreover, the set of closed geodesic polygons is denoted by $\Pcl(M,\T) :=
\bigcup_{x\in M}\P_x^x(M,\T)$.

With the partition $\T$ being given, the correspondence $\gamma \leftrightarrow
(x_0,\ldots,x_r)$ identifies $\P(M,\T)$ with an open dense subset of
$\underbrace{M \times \cdots \times M}_{r+1\,\, \mathrm{factors}}$.
The complement of this subset is a zero-set.
The Riemannian product volume measure on $M \times \cdots \times M$ induces a
measure on $\P(M,\T)$ which we denote by $\Dg$.
Similarly, $\P_x(M,\T)$, $\P^y(M,\T)$, and $\Pcl(M,\T)$ can be identified with 
dense open subsets of $\underbrace{M \times \cdots \times M}_{r\,\,
  \mathrm{factors}}$.
The corresponding measures on these spaces will also be denoted by $\Dg$.

Recall that the {\em energy} of a piecewise smooth curve $\gamma:[0,t] \to M$
is defined as 
$$
\E(\gamma) := \frac12 \int_0^t |\dot\gamma(s)|^2 ds.
$$
{\em Parallel transport} in $E$ along $\gamma$ with respect to the connection
$\nabla$ will be denoted by $\tau(\gamma,\nabla) : E_{\gamma(0)} \to
E_{\gamma(t)}$.
More generally, for $s,s' \in [0,t]$ we have parallel transport
$\tau(\gamma,\nabla)_s^{s'} : E_{\gamma(s)} \to E_{\gamma(s')}$.
Parallel transport is a linear isometry.
In particular, its operator norm equals $1$.
We have $\tau(\gamma,\nabla)_{s'}^{s''}\circ\tau(\gamma,\nabla)_s^{s'} =
\tau(\gamma,\nabla)_s^{s''}$ and $(\tau(\gamma,\nabla)_s^{s'})^{-1} =
\tau(\gamma,\nabla)_{s'}^{s}$.

Finally, we put
\begin{equation}
Z(\T,m) := \prod_{j=1}^r(4\pi t_j)^{m/2} .
\label{eq:ZTerm}  
\end{equation}
Now we can state the main result.

\begin{thm}\label{thm:PathIntegral1}
Let $M$ be an $m$-dimensional closed Riemannian manifold, let $E$ be a vector
bundle over $M$ with a metric and a compatible connection $\nabla$.
Let $H = \nabla^*\nabla + V$ be a self-adjoint generalized Laplace
operator acting on sections in $E$.
Let $t>0$.

Then for any sequence of partitions $\T_n =(t_1^n, \ldots,t_{r_n}^n)$ with
$|\T_n| \to 0$ and $\Laenge(\T_n) \to t$ as $n \to \infty$ and for any $u\in
C^0(M,E)$ 
\begin{equation*}
\frac{1}{Z(\T_n,m)} \,
\int_{\P_x(M,\T_n)} 
\exp\left(-\frac12\E(\gamma)+\int_0^{\Laenge(\T_n)}\frac13 \scal(\gamma(s))\, 
ds\right)\cdot
\tau(\gamma,\nabla)_{\Laenge(\T_n)}^{0}\times
\end{equation*}
\begin{equation*}
\times
\prod_{j=1}^{r_n}
\exp\left(-\int_{\sigma_{j-1}(\T_n)}^{\sigma_j(\T_n)}
  \tau(\gamma,\nabla)_s^{\Laenge(\T_n)}\cdot
  V(\gamma(s)) 
  \cdot\tau(\gamma,\nabla)_{\Laenge(\T_n)}^{s}\,ds\right)\cdot 
u(\gamma(\Laenge(\T_n)))\,\Dg
\end{equation*}
\begin{equation*}
\quad\quad\xrightarrow{n \to \infty} \quad\quad e^{-tH}u(x)
\end{equation*}
converges uniformly in $x$.
\end{thm}

The proof will be provided in the next section.

\begin{cor}\label{cor:PathIntegral1}
If in addition to the assumptions in Theorem~\ref{thm:PathIntegral1} the
potential $V(x)$ is a scalar multiple of the identity for each $x\in M$, then 
\begin{equation*}
\frac{1}{Z(\T_n,m)} \, \int_{\P_x(M,\T_n)} 
\exp\left(-\frac12\E(\gamma) +\int_{0}^{L(\T_n)} \left(\frac13 \scal(\gamma(s))
- V(\gamma(s)) \right)ds\right)\times
\end{equation*}
\begin{equation*}
\times\tau(\gamma,\nabla)_{\Laenge(\T_n)}^{0}\cdot u(\gamma(\Laenge(\T_n)))\,\Dg
\end{equation*}
\begin{equation*}
\quad\quad\quad\quad
\xrightarrow{n \to \infty} \quad\quad e^{-tH}u(x)
\end{equation*}
converges uniformly in $x$.
\end{cor}

\begin{proof}[Proof of Corollary]
If $V$ is scalar, then all operators in the following integrals commute and we
have 
\begin{eqnarray*}
\lefteqn{\hspace{-1.5cm}\prod_{j=1}^{r_n}
\exp\int_{\sigma_{j-1}(\T_n)}^{\sigma_j(\T_n)} \left(
-\tau(\gamma,\nabla)_s^{\Laenge(\T_n)}\cdot V(\gamma(s))
  \cdot\tau(\gamma,\nabla)_{\Laenge(\T_n)}^{s}\right)ds}\\
&=&
\prod_{j=1}^{r_n}
\exp\int_{\sigma_{j-1}(\T_n)}^{\sigma_j(\T_n)} \left(
- V(\gamma(s))\right)ds\\
&=&
\exp\sum_{j=1}^{r_n}
\int_{\sigma_{j-1}(\T_n)}^{\sigma_j(\T_n)} \left(
- V(\gamma(s))\right)ds\\
&=&
\exp
\int_{0}^{L(\T_n)} \left(
- V(\gamma(s))\right)ds .
\end{eqnarray*}
Theorem~\ref{thm:PathIntegral1} implies the corollary.
\end{proof}

\section{Proof of the main theorem}

This section is devoted to the proof of Theorem~\ref{thm:PathIntegral1}.
For $t>0$ the operator $e^{-tH}$ has an integral kernel $k_t$, i.~e.,
$$
e^{-tH}u(x) = \int_M k_t(x,y)\,u(y)\,dy .
$$
This integral kernel $(t,x,y)\mapsto k_t(x,y)$ is smooth on $(0,\infty)\times
M\times M$. 
Hence the solution $U(t,x)=e^{-tH}u(x)$ of the heat equation
$\frac{\partial U}{\partial t}+HU=0$ with initial condition $U(0,x)=u(x)$ is
smooth on $(0,\infty)\times M$ for any $u\in L^2(M,E)$.
If $u$ is continuous, $u\in C^0(M,E)$, then
the solution $U(t,x)$ is continuous at $t=0$ as well.
Hence $t\mapsto e^{-tH}$ yields a strongly continuous semigroup in the Banach
space $C^0(M,E)$ of continuous sections in $E$.

For the proof we proceed by successively replacing the heat semigroup by
Chernoff equivalent proper families.
This will be done by suitably modifying the integral kernels.
We will apply the results from Section~\ref{sec:cherno} for proper families $S$
with $DS=-H$ on the Banach space $X=C^0(M,E)$. 

Our analysis of the solutions to the heat equation will be based on the
precise understanding of the short time asymptotics of the heat kernel.
To formulate the result let $\chi :[0,\infty)\to[0,1]$ be a monotonic smooth
function being $1$ near $0$ and with support in $[0,\injrad(M)^2/4)$.
Here $\injrad(M)$ denotes the \emph{injectivity radius} of $M$.
We have (see \cite[Thm.~2.30, p.~87]{BGV})

\begin{thm}[Heat kernel asymptotics]\label{thm:HeatAsymptotics}
There exist smooth sections $\Phi_j\in C^\infty(M\times M,E\boxtimes E^*)$
such that the kernel defined by
\begin{equation}\label{asymptkern}
{k}_t^{(n)}(x,y)
:=
\left({4\pi t}\right)^{-m/2}\cdot
\exp\left(-\frac{d(x,y)^2}{4t}\right)\cdot
\chi\left(d(x,y)^2\right)\cdot
\sum_{j=0}^n t^j\cdot\Phi_j(x,y)
\end{equation}
is asymptotic to the heat kernel $k_t(x,y)$ in the sense that for all
$n>\tfrac{m}{2}$ one has 
\begin{equation}\label{asymptquantitativ}
\left\|k_t-{k}_t^{(n)}\right\|_{C^0(M\times M)}
= 
\OOOO\left(t^{n-m/2}\right)\qquad\mbox{ as }\quad t\searrow 0.
\end{equation}
\end{thm}
\vspace{-0.8cm}
\hfill$\Box$
\vspace{0.8cm}

Here $E\boxtimes E^*$ denotes the exterior tensor product whose fiber over
$(x,y) \in M\times M$ is given by $(E\boxtimes E^*)_{(x,y)} = E_x \otimes
E^*_y = \Hom(E_y,E_x)$.
Because of the cut-off factor $\chi\left(d(x,y)^2\right)$ the section $(x,y)
\mapsto \Phi_j(x,y)$ needs to be specified only for $d(x,y) < \injrad(M)$.
In particular, in this case $x$ and $y$ are not cut-points of each other, thus
there is a unique (up to reparametrization) shortest geodesic joining $x$ and
$y$.

For any $y\in M$ the Riemannian exponential map $\exp_y : T_yM \to M$ is
smooth and onto.
It maps a star-shaped open neighborhood of $0\in T_yM$ diffeomorphically onto
the complement of the cut-locus $\CC_y$ of $y$.
For $x\in M\setminus \CC_y$ we put 
$$
\mu(x,y) := \det(d\exp_y(\exp_y^{-1}(x))).
$$
The function $\mu$ measures the volume distortion of the Riemannian
exponential map.
In Riemannian normal coordinates about $y$ we have $\mu(x,y) =
\sqrt{\det(g_{ij})}$. 
The function $\mu$ is smooth and positive on its domain $\{(x,y)\in M\times
M\, |\, x \not\in \CC_y\}$.
Since $\det(d\exp)$ is smooth on all of $TM$ and $M$ is compact $\mu$ is also
bounded from above.
Moreover, we have $\mu(x,x)=1$.

If $d(x,y)<\injrad(M)$ let $\gamma_{x,y}$ be a shortest geodesic joining $x$
and $y$.
The leading term $\Phi_0(x,y)$ is given by 
$$
\Phi_0(x,y) = \mu(x,y)^{-1/2} \cdot \tau(\gamma_{x,y},\nabla).
$$
In particular, $\Phi_0(x,x)=\id$.

\subsection{First Modification}
For $t>0$ denote the integral operator with integral kernel $k^{(n)}_t$ by
$K^{(n)}_t$, i.~e.,
$$
K^{(n)}_t\,u(x) := \int_Mk^{(n)}_t(x,y)\,u(y)\, dy .
$$
For $t=0$ put $K^{(n)}_0:=\id$.

\begin{lem}\label{lem:Cutoff1}
For $n > 1+\frac{m}{2}$ the family $[0,\infty) \to
\B(C^0(M,E))$, $t\mapsto K^{(n)}_t$, is proper and Chernoff equivalent to
the semigroup $(e^{-tH})_{t\geq0}$.
\end{lem}

\begin{proof}
Let $u\in C^0(M,E)$.
Then 
\begin{eqnarray*}
\|e^{-tH}u - K^{(n)}_t u\|_{C^0(M)}
&=&
\left\|\,\int_M\left(k_t(\cdot,y)-k^{(n)}_t(\cdot,y)\right)
u(y)\,dy\right\|_{C^0(M)}\\
&\leq&
\vol(M)\cdot \left\|k_t-{k}_t^{(n)}\right\|_{C^0(M\times M)} 
\cdot \|u\|_{C^0(M)}
\end{eqnarray*}
holds for $t>0$.
From Theorem~\ref{thm:HeatAsymptotics} we conclude
$$
\|e^{-tH} - K^{(n)}_t\|_{C^0,C^0} 
\quad \leq\quad  
\vol(M)\cdot \left\|k_t-{k}_t^{(n)}\right\|_{C^0(M\times M)}\\
\quad =\quad  \OOOO\left(t^{n-m/2}\right) .
$$
Since $n-m/2>1$ we have 
\begin{equation}
\|e^{-tH} - K^{(n)}_t\|_{C^0,C^0}=\oooo(t)
\label{eq:etHKm}
\end{equation}
The statement now follows from Lemma~\ref{lem:ChernoffKriterium}.
\end{proof}

Put 
$$
e_t(x,y) := \left({4\pi t}\right)^{-m/2}\cdot
\exp\left(-\frac{d(x,y)^2}{4t}\right)\cdot
\chi\left(d(x,y)^2\right) .
$$

The following lemma will allow us to perform the next modifications. 

\begin{lem}[Workhorse lemma]\label{lem:KernKriterium}
Let $p_t(x,y)$ and $q_t(x,y) \in E_x \otimes E_y^*$ depend continuously on
$(t,x,y) \in (0,\infty) \times M \times M$.
Denote the corresponding integral operators by $P_t$ and $Q_t$ respectively.
We assume that $P_0:= \id$ extends $(P_t)_{t>0}$ to a proper family
$(P_t)_{t\geq0}$ in $C^0(M,E)$.
Suppose there exist constants $C,\alpha,\beta \geq 0$, $\beta+\alpha/2 > 1$,
such that
$$
|p_t(x,y) - q_t(x,y)| 
\quad\leq \quad
C \cdot e_t(x,y) \cdot d(x,y)^\alpha \cdot t^\beta
$$
for all $0<t\leq t_0$ and for all $x,y \in M$ where $t_0$ is some
positive constant.

Then $Q_0:= \id$ extends $(Q_t)_{t>0}$ to a proper family $(Q_t)_{t\geq0}$ in
$C^0(M,E)$ which is Chernoff equivalent to $(P_t)_{t\geq0}$.
\end{lem}

\begin{proof}
By Lemma~\ref{lem:ChernoffKriterium} it is sufficient to show  
\begin{equation}
\|P_t-Q_t\|_{C^0,C^0} = \oooo(t)
\label{eq:PtQt}
\end{equation}
as $t\searrow 0$.
We choose a constant $C_1>0$ such 
that $\tau^\alpha \leq C_1 \cdot
\exp(\tau^2)$ for all $\tau \in [0,\infty)$.
With $\tau = {d(x,y)}/{\sqrt{8t}}$ this yields
\begin{equation}
d(x,y)^\alpha 
\quad\leq\quad
C_1 \cdot (8t)^{\alpha/2} \cdot \exp\left(\frac{d(x,y)^2}{8t}\right) .
\label{est:dxygegent}
\end{equation}
Hence
\begin{eqnarray*}
|p_t(x,y) - q_t(x,y)|
&\leq&
C \cdot e_t(x,y)\cdot d(x,y)^\alpha \cdot t^\beta\\
&\stackrel{(\ref{est:dxygegent})}{\leq}&
C_2 \cdot e_t(x,y)\cdot t^{\beta+\alpha/2} \cdot
\exp\left(\frac{d(x,y)^2}{8t}\right)  \\
&=&
C_3 \cdot e_{2t}(x,y)\cdot t^{\beta+\alpha/2} .
\end{eqnarray*}
We fix $x\in M$ and consider Riemannian normal coordinates given by the
Riemannian exponential map $\exp_x : T_xM \to M$.
Due to the cut-off term in $e_t(x,y)$ the integrand is supported in the region
where $\exp_x$ is a diffeomorphism.
For any $u\in C^0(M,E)$ we have 
\begin{eqnarray*}
\lefteqn{|(P_t-Q_t)u(x)|}\\
&\leq&
C_3 \cdot t^{\beta+\alpha/2} \cdot \int_M e_{2t}(x,y) \cdot |u(y)| \, dy \\
&=&
C_3 \cdot t^{\beta+\alpha/2} \cdot \int_{T_xM}
\frac{e^\frac{-|\xi|^2}{8t}}{(8\pi t)^{m/2}} \cdot
|u(\exp_x(\xi))| \cdot \chi(|\xi|^2) \cdot |\det (d\exp_x(\xi))| \,  d\xi \\
&\leq&
C_4 \cdot t^{\beta+\alpha/2} \cdot \|u\|_{C^0(M)} \cdot  
\int_{T_xM}\frac{e^\frac{-|\xi|^2}{8t}}{(8\pi t)^{m/2}}\,  d\xi \\
&=&
C_4 \cdot t^{\beta+\alpha/2} \cdot \|u\|_{C^0(M)} ,
\end{eqnarray*}
the constant $C_4$ being independent of $x$.
Since $\beta+\alpha/2>1$ this shows (\ref{eq:PtQt}) and proves the lemma.
\end{proof}

\subsection{Second modification}
We replace $K_t^{(n)}$ by $K_t^{(1)}$.

\begin{lem}\label{lem:Cutoff2}
The family $[0,\infty) \to \B(C^0(M,E))$, $t\mapsto K^{(1)}_t$, is proper and
Chernoff equivalent to the semigroup $(e^{-tH})_{t\geq0}$.
\end{lem}

\begin{proof}
By Lemma~\ref{lem:Cutoff1} it is sufficient to show that $(K_t^{(1)})$ is a
proper family which is Chernoff equivalent to $(K_t^{(n)})$.  
From
$$ 
k^{(n)}_t(x,y) - k^{(1)}_t(x,y) 
= 
e_t(x,y)\cdot \sum_{j=2}^n t^j\cdot\Phi_j(x,y)
$$
we see that Lemma~\ref{lem:KernKriterium} applies with $P_t=K^{(n)}_t$,
$Q_t=K^{(1)}_t$, $\alpha=0$ and $\beta=2$. 
\end{proof}

\subsection{Third modification}
The integral kernel of $K^{(1)}_t$ is given by
\begin{eqnarray*}
k^{(1)}_t(x,y) 
&=&
e_t(x,y)\cdot \left(\Phi_0(x,y)+t\cdot\Phi_1(x,y)\right)\\
&=&
e_t(x,y)\cdot \Phi_0(x,y)\cdot
\left(\id+t\cdot\Phi_0(x,y)^{-1}\cdot\Phi_1(x,y)\right)\\
&=&
e_t(x,y)\cdot \Phi_0(x,y)\cdot
\left(\id+t\cdot A(x,y)\right)
\end{eqnarray*}
where we have put $A(x,y):=\Phi_0(x,y)^{-1}\cdot\Phi_1(x,y)\in\Hom(E_y,E_y)$.
Recall that $\Phi_0(x,y)$ is $\mu(x,y)^{-1/2}$ times parallel transport along the
shortest geodesic from $y$ to $x$ and hence is invertible.
Also remember that due to the cut-off factor $\chi(d(x,y)^2)$ in $e_t(x,y)$
there is no need to specify $\Phi_0(x,y)$ if $y$ lies in the cut-locus of $x$
in which case the shortest geodesic from $y$ to $x$ may not be unique and
$\mu(x,y)$ is not defined.
We define a new smooth integral kernel by
$$
w_t(x,y) := e_t(x,y)\cdot \Phi_0(x,y)\cdot\exp\left(t\cdot A(x,y)\right) .
$$
Here $\exp\left(t\cdot A(x,y)\right)$ is the usual exponential of
endomorphisms of $E_y$ defined by the power series $\exp\left(t\cdot
  A(x,y)\right) = \sum_{j=0}^\infty \frac{t^j}{j!}A(x,y)^j$.
The corresponding integral operator is denoted by $W_t$,
$$
W_tu(x) = \int_Mw_t(x,y)\,u(y)\,dy .
$$

\begin{lem}\label{lem:Cutoff3}
The family $[0,\infty) \to \B(C^0(M,E))$, $t\mapsto W_t$, is proper and
Chernoff equivalent to the semigroup $(e^{-tH})_{t\geq0}$.
\end{lem}

\begin{proof}
This follows from Lemmas~\ref{lem:Cutoff2} and \ref{lem:KernKriterium} with
$P_t=K^{(1)}_t$, $Q_t=W_t$, $\alpha=0$, and $\beta=2$ because
$$
w_t(x,y) - k^{(1)}_t(x,y) 
\quad=\quad 
e_t(x,y) \cdot \Phi_0(x,y) \cdot \OOOO(t^2).
$$
\end{proof}

\subsection{Fourth modification}
For $x,y\in M$ with $d(x,y)< \injrad(M)$ and $a<b$ we denote
$\gabxy$ the unique minimal geodesic parametrized on $[a,b]$ with
$\gabxy(a)=x$ and $\gabxy(b)=y$.
We put for $t>0$ and $x,y\in M$
\begin{eqnarray*}
v_t(x,y) 
&:=& 
e_t(x,y)\cdot \Phi_0(x,y) \cdot
\exp\left(\int_0^t \tau(\gamma,\nabla)_s^t\cdot A(\gamma(s),\gamma(s))
\cdot\tau(\gamma,\nabla)_t^s\,ds\right)\\
&=&
e_t(x,y)\cdot \mu(x,y)^{-1/2}\cdot\tau(\gamma,\nabla)_t^0\times\\
&&
\times\exp\left(\int_0^t \tau(\gamma,\nabla)_s^t\cdot A(\gamma(s),\gamma(s))
\cdot\tau(\gamma,\nabla)_t^s\,ds\right)
\end{eqnarray*}
where $\gamma=\gotxy$.
Correspondingly, we set 
$$
V_tu(x) := \int_Mv_t(x,y)\, u(y)\, dy.
$$

\begin{lem}\label{lem:Cutoff4}
Putting $V_0:=\id$ we obtain a proper family $(V_t)_{t\ge0}$ being Chernoff
equivalent to the semigroup $(e^{-tH})_{t\geq0}$. 
\end{lem}

\begin{proof}
For $x$ and $y$ in $M$ we have that the distance of $(x,y)$ and
$(\gotxy(s),\gotxy(s))$ in $M \times M$ is bounded from above by $C_1\cdot
d(x,y)$ for all $s\in[0,t]$.
Thus, abreviating $\gamma=\gotxy$,
$$
\left| A(x,y) - 
\tau(\gamma,\nabla)_s^t \cdot A(\gamma(s),\gamma(s)) \cdot
\tau(\gamma,\nabla)_t^s \right| \leq C_2 \cdot d(x,y), 
$$
hence
\begin{eqnarray*}
\lefteqn{\left| t\cdot A(x,y) - \int_0^t
  \tau(\gamma,\nabla)_s^t \cdot A(\gamma(s),\gamma(s)) \cdot
  \tau(\gamma,\nabla)_t^s\, ds  \right|}\\ 
&\leq&
\int_0^t \left| A(x,y) - \tau(\gamma,\nabla)_s^t  \cdot A(\gamma(s),\gamma(s))
  \cdot \tau(\gamma,\nabla)_t^s  \right|\, ds\\ 
&\leq& 
C_2 \cdot d(x,y)\cdot t
\end{eqnarray*}
and therefore
$$
\left| \exp\left(t\cdot A(x,y)\right) 
- \exp\left(\int_0^t \tau(\gamma,\nabla)_s^t \cdot A(\gamma(s),\gamma(s))
\cdot \tau(\gamma,\nabla)_t^s\, ds\right)  \right| 
\leq
C_3 \cdot d(x,y)\cdot t.
$$
Here we applied a local Lipschitz bound for $\exp$ which is justified since
$t$ and $A$ are bounded.
We get
\begin{eqnarray*}
\lefteqn{|w_t(x,y)-v_t(x,y)|}\\
&\leq&
e_t(x,y)\cdot |\Phi_0(x,y)|\times \\
&&
\times
\left| \exp\left(t\cdot A(x,y)\right) 
- \exp\left(\int_0^t \tau(\gamma,\nabla)_s^t \cdot A(\gamma(s),\gamma(s))
\cdot \tau(\gamma,\nabla)_t^s\, ds\right)  \right| \\ 
&\leq&
e_t(x,y)\cdot C_4 \cdot d(x,y)\cdot t .
\end{eqnarray*}
The assertion follows from Lemma~\ref{lem:Cutoff3} and
Lemma~\ref{lem:KernKriterium} with $P_t=W_t$, $Q_t=V_t$, and $\alpha=\beta=1$.
\end{proof}

The advantage of the integral kernel of $V_t$ compared to the one of $W_t$
lies in the fact that we need to evaluate $A$ only along the diagonal.
It is well-known that 
\begin{equation}
A(x,x) = \Phi_0(x,x)^{-1}\cdot \Phi_1(x,x) 
= \Phi_1(x,x) = \frac16 \scal(x) \cdot\id - V(x),
\label{eq:A}
\end{equation}
compare the computations in~\cite[p.~103ff]{Roe}.
Here $\scal$ denotes the scalar curvature of $M$ and $V$ the potential of $H$.
Thus the integral kernel of $V_t$ is given by 
\begin{eqnarray}
v_t(x,y) 
&=&
e_t(x,y)\cdot \mu(x,y)^{-1/2}\cdot\tau(\gamma,\nabla)_t^0\times\nonumber\\
&&
\times\exp\int_0^t \left(\frac16 \scal(\gamma(s))\cdot\id -
  \tau(\gamma,\nabla)_s^t\cdot V(\gamma(s))
  \cdot\tau(\gamma,\nabla)_t^s\right) ds
\label{eq:KernVt}
\end{eqnarray}
where $\gamma=\gamma_{x,y}^{0,t}$.

\subsection{Fifth modification}
Next we replace the term $\mu(x,y)^{-1/2}$ by an additional scalar curvature
term. 
We put 
\begin{eqnarray*}
\lefteqn{\hvt(x,y)}\\ 
&:=&
e_t(x,y)\cdot\tau(\gamma,\nabla)_t^0\cdot
\exp\int_0^t \left(\frac13 \scal(\gamma(s))\cdot\id -
  \tau(\gamma,\nabla)_s^t\cdot V(\gamma(s))
  \cdot\tau(\gamma,\nabla)_t^s\right) ds 
\end{eqnarray*}
and
$$
\hVt u (x) := \int_M \hvt(x,y)\cdot u(y) \, dy .
$$
\begin{lem}\label{lem:Cutoff5}
Putting $\widehat V_0:=\id$ we obtain a proper family $(\hVt)_{t\ge0}$ being
Chernoff equivalent to the semigroup $(e^{-tH})_{t\geq0}$. 
\end{lem}

\begin{proof}
We can compute in Riemannian normal coordinates about $y$ because the cut-off
term in $e_t(x,y)$ ensures that only points $x$ in the domain of this
coordinate system matter.
Putting $\xi := \exp_y^{-1}(x)$ it is well-known that the determinant of the
metric has the Taylor expansion
$$
\det g = 1 - \frac13 \cdot \ric_y(\xi,\xi) + \OOOO(|\xi|^3)
$$
where the constant in the $\OOOO(|\xi|^3)$ term can be chosen independently of
$y$, compare \cite[Cor.~2.3, p.~84]{C}.
Hence 
$$
\mu(x,y)^{-1/2} = (\det g)^{-1/4} = 1 + \frac{1}{12} \cdot \ric_y(\xi,\xi) +
\OOOO(|\xi|^3) .
$$
Since $|\xi|=d(x,y)$ it follows from Lemma~\ref{lem:KernKriterium} with
$\alpha=3$ and $\beta=0$ that we can replace $\mu(x,y)^{-1/2}$ in the kernel
$v_t(x,y)$ by $1+ \frac{1}{12} \cdot \ric_y(\exp_y^{-1}(x),\exp_y^{-1}(x))$ to
obtain a new proper family of integral operators being Chernoff equivalent to
$(V_t)_t$ (and hence to $e^{-tH}$).

Next we apply Lemma~\ref{lem:GaussEst} in normal coordinates about $y$ with
$$
B\quad=\quad\frac{1}{12}\cdot \ric_y
$$ 
and 
\begin{eqnarray*}
\lefteqn{f(t,\xi) \quad=\quad 
\chi(|\xi|^2) \cdot\tau(\gamma,\nabla)_t^0 \times}\\
&&
\times \exp\int_0^t \left(\frac16 \scal(\exp_y(s\xi))\cdot \id 
- \tau(\gamma,\nabla)_s^t
  V(\exp_y(s\xi))\tau(\gamma,\nabla)_t^s\right)ds \cdot u(\exp_y(\xi)).
\end{eqnarray*}
In order to apply Lemma~\ref{lem:ChernoffKriterium1} we only need to insert
$u$ of the form $u=e^{-aH}v$ for $a>0$.
Since $e^{-aH}$ is smoothing all such $u$ are smooth.
Thus in this case $f$ is smooth and has compact support in $\xi$.
Therefore Lemma~\ref{lem:GaussEst} applies and
Lemma~\ref{lem:ChernoffKriterium1} tells us that we may replace
$B(\xi,\xi)=\frac{1}{12}\cdot \ric_y(\xi,\xi)$ by $2t \cdot
\tr(\frac{1}{12}\cdot \ric_y) = \frac{t}{6}\scal(y)$. 
We obtain a new proper family of integral operators $\tVt$, Chernoff
equivalent to $V_t$, with integral kernel 
\begin{eqnarray*}
\tvt(x,y) 
&=&
e_t(x,y)\cdot \left(1+\frac{t}{6}\scal(y)\right)\cdot
\tau(\gamma,\nabla)_t^0\times\nonumber\\ 
&&
\times\exp\int_0^t \left(\frac16 \scal(\gamma(s))\cdot\id -
  \tau(\gamma,\nabla)_s^t\cdot V(\gamma(s))
  \cdot\tau(\gamma,\nabla)_t^s\right) ds .
\end{eqnarray*}
Since
$$
\left| \frac{t}{6}\scal(y) - \int_0^t \frac16 \scal(\gotxy(s))\,ds\right|
\quad\leq\quad 
C \cdot t \cdot d(x,y)
$$
Lemma~\ref{lem:KernKriterium} applies with $\alpha=\beta=1$ and shows that we
may replace $\frac{t}{6}\scal(y)$ by $\int_0^t \frac16 \scal(\gotxy(s))\,ds$
to get another proper family of Chernoff equivalent integral operators with
integral kernel
\begin{eqnarray*}
\bvt(x,y) 
&=&
e_t(x,y)\cdot \left(1+ \int_0^t \frac16 \scal(\gotxy(s))\,ds\right)\cdot
\tau(\gamma,\nabla)_t^0\times\nonumber\\ 
&&
\times\exp\int_0^t \left(\frac16 \scal(\gamma(s))\cdot\id -
  \tau(\gamma,\nabla)_s^t\cdot V(\gamma(s))
  \cdot\tau(\gamma,\nabla)_t^s\right) ds .
\end{eqnarray*}
Finally, since
$$
\left(1+\int_0^t \frac16 \scal(\gotxy(s))\,ds\right) -
\exp\left(\int_0^t \frac16 \scal(\gotxy(s))\,ds\right)
\quad=\quad 
\OOOO(t^2)
$$
Lemma~\ref{lem:KernKriterium} with $\alpha=0$ and $\beta=2$ yields the
Chernoff equivalent proper family of integral operators with integral kernel
\begin{eqnarray*}
\lefteqn{\hvt(x,y)}\\ 
&=&
e_t(x,y)\cdot \exp\left(\int_0^t \frac16 \scal(\gamma(s))\,ds\right)\cdot
\tau(\gamma,\nabla)_t^0\times \\ 
&&
\times\exp\int_0^t \left(\frac16 \scal(\gamma(s))\cdot\id -
  \tau(\gamma,\nabla)_s^t\cdot V(\gamma(s))
  \cdot\tau(\gamma,\nabla)_t^s\right) ds \\
&=&
e_t(x,y)\cdot\tau(\gamma,\nabla)_t^0\cdot
\exp\int_0^t \left(\frac13 \scal(\gamma(s))\cdot\id -
  \tau(\gamma,\nabla)_s^t\cdot V(\gamma(s))
  \cdot\tau(\gamma,\nabla)_t^s\right) ds .
\end{eqnarray*}
\end{proof}

\begin{rem}\label{rem:scal}
Lemma~\ref{lem:Cutoff5} can be generalized as follows.
Instead of $\mu(x,y)^{-1/2}$ we fix $\Lambda\in\R$ and we consider 
$$
\mu(x,y)^{-\Lambda/2} = (\det g)^{-\Lambda/4} 
= 1 + \frac{\Lambda}{12} \cdot \ric_y(\xi,\xi) + \OOOO(|\xi|^3) .
$$
Now the same proof yields a proper family with integral kernel
$$
e_t(x,y)\cdot\mu(x,y)^{\frac{\Lambda-1}{2}}\cdot\tau(\gamma,\nabla)_t^0\times
$$
$$
\times\exp\int_0^t \left(\frac{\Lambda+1}{6} \scal(\gamma(s))\cdot\id -
  \tau(\gamma,\nabla)_s^t\cdot V(\gamma(s))
  \cdot\tau(\gamma,\nabla)_t^s\right) ds 
$$
being Chernoff equivalent to the semigroup $(e^{-tH})_{t\geq 0}$.
\end{rem}

\subsection{Sixth modification}
Finally, we remove the cut-off term.
Put for $t>0$ and $x,y\in M$, $y\not\in \CC_x$,
\begin{eqnarray*}
\hwt(x,y)
&:=& 
\left({4\pi t}\right)^{-m/2}\cdot
\exp\left(-\frac{d(x,y)^2}{4t}\right)\cdot\tau(\gamma,\nabla)_t^0\times\\
&&
\times\exp\int_0^t \left(\frac13 \scal(\gamma(s))\cdot\id -
  \tau(\gamma,\nabla)_s^t\cdot V(\gamma(s))
  \cdot\tau(\gamma,\nabla)_t^s\right) ds 
\end{eqnarray*}
where $\gamma=\gotxy$ and 
$$
\hWt u(x) := \int_{M\setminus\CC_x}\hwt(x,y)\, u(y)\, dy.
$$

\begin{lem}\label{lem:Cutoff6}
Putting $\widehat W_0:=\id$ we obtain a proper family $(\hWt)_{t\ge0}$ being
Chernoff equivalent to the semigroup $(e^{-tH})_{t\geq0}$.
\end{lem}

\begin{proof}
We show $\|\hWt - \hVt\| = \OOOO(t^2)$.
The result then follows from Lemmas~\ref{lem:ChernoffKriterium} and
\ref{lem:Cutoff5}.
We have
$$
|\hwt(x,y) - \hvt(x,y)| 
\quad\leq\quad
C_1\cdot\left({4\pi t}\right)^{-m/2}\cdot
\exp\left(-\frac{d(x,y)^2}{4t}\right)\cdot (1-\chi(d(x,y)^2)) .
$$
Since $1-\chi(d(x,y)^2)=0$ whenever $d(x,y)^2\leq \eps$ for some suitable
$\eps > 0$ we get for small $t>0$
\begin{eqnarray*}
|\hwt(x,y) - \hvt(x,y)|
&\leq&
C_1\cdot\left({4\pi t}\right)^{-m/2}\cdot
\exp\left(-\frac{\eps}{4t}\right)\\
&\leq&
C_2\cdot\left({4\pi t}\right)^{-m/2}\cdot t^{2+m/2}\\
&=&
C_3 \cdot t^2.
\end{eqnarray*}
The lemma follows.
\end{proof}

\subsection{Conclusion of proof}
As a consequence of Chernoff's Theorem~\ref{thm:Chernoff} we 
know that for any sequence of partitions
$\T_n=(t_1^n,\ldots,t_{r_n}^n)$ satisfying
\begin{equation}
\lim_{n\to\infty} \Laenge(\T_n) = t
\quad\quad\mbox{ and }\quad
\lim_{n\to\infty} |\T_n| = 0
\label{eq:Partition}
\end{equation}
and any $u\in C^0(M,E)$ we have
\begin{equation}
\hW_{t_{1}^n}\cdots \hW_{t_{r_n}^n}u \rightarrow e^{-tH}u 
\label{eq:CKonv}
\end{equation}
uniformly as $n\to\infty$.
We will study the operators $\hW_{t_{1}^n}\cdots \hW_{t_{r_n}^n}$ in
more detail.
Consider the Riemannian manifold $\M^r := \{(x_1,\ldots,x_r) \in
M\times\cdots\times M\, |\, x_{j+1} \not\in \CC_{x_j}\}$ equipped with the
product metric. 
For fixed $x_0\in M$ we may write, using the abbreviation
$X=(x_1,\ldots,x_{r_n})$ and $\gamma_j=\gamma_{x_{j-1},x_j}^{0,t_j^n}$, 
\begin{eqnarray}
\lefteqn{\hW_{t_{1}^n}\cdots \hW_{t_{r_n}^n}u(x_0)}\nonumber \\
&=&
\int_{\M^{r_n}} \hw_{t_{1}^n}(x_0,x_1)\cdots 
\hw_{t_{r_n}^n}(x_{r_n-1},x_{r_n})\,u(x_{r_n})\,dX\nonumber \\
&=&
(4\pi)^{-r_nm/2}\cdot\prod_{j=1}^{r_n}(t_j^n)^{-m/2} \cdot 
\int_{\M^{r_n}} 
\exp\left(-\sum_{j=1}^{r_n}\frac{d(x_{j-1},x_j)^2}{4t_j^n}\right) \cdot
\prod_{j=1}^{r_n}\Big[\tau(\gamma_j,\nabla)_0^{t_j^n}\times\nonumber\\
&&
\times\exp\int_0^{t_j^n}\left(
\frac13\scal(\gamma_j(s))-\tau(\gamma_j,\nabla)_s^{t_j^n}\cdot 
  V(\gamma_j(s))
  \cdot\tau(\gamma_j,\nabla)_{t_j^n}^{s}\right)ds\Big]
\cdot u(x_{r_n})\,dX .
\label{eq:Expand1}
\end{eqnarray}

For any partition $\T=(t_1,\ldots,t_r)$ we write as before $\sigma_j(\T) =
t_1+\ldots +t_j$. 
To $x_0\in M$ and $X=(x_1,\ldots,x_r)\in\M^r$
we associate the geodesic polygon $\gamma_{x_0,X}^\T :
[0,\Laenge(\T)] \to M$ which, when restricted to the subinterval
$[\sigma_{j-1}(\T),\sigma_{j}(\T)]$, is the shortest curve from
$x_{j-1}$ to $x_j$, parametrized proportionally to arclength.
In other words, we have
$\gamma_{x_0,X}^\T|_{[\sigma_{j-1}(\T),\sigma_{j}(\T)]} =
\gamma_{x_{j-1},x_j}^{\sigma_{j-1}(\T),\sigma_j(\T)}$. 
We get 
\begin{eqnarray}
\lefteqn{\prod_{j=1}^{r}\Big[\tau(\gamma_j,\nabla)_0^{t_j}\cdot
\exp\int_0^{t_j}
\left(\frac13\scal(\gamma_j(s))-\tau(\gamma_j,\nabla)_s^{t_j}\cdot 
  V(\gamma_j(s))
  \cdot\tau(\gamma_j,\nabla)_{t_j}^{s}\right)ds\Big]}\nonumber\\
&=&
\exp\left(\int_0^{\Laenge(\T)}\frac13 \scal(\gamma_{x_0,X}^\T(s))\, ds\right)
\cdot\tau(\gamma_{x_0,X}^\T,\nabla)_{\Laenge(\T)}^{0}\times
\nonumber\\
&&
\times\prod_{j=1}^{r}
\exp\left(-\int_{\sigma_{j-1}(\T)}^{\sigma_j(\T)} 
\tau(\gamma_{x_0,X}^\T,\nabla)_s^{\Laenge(\T)}\cdot V(\gamma_{x_0,X}^\T(s))
  \cdot\tau(\gamma_{x_0,X}^\T,\nabla)_{\Laenge(\T)}^{s}\,ds\right) .
\label{eq:IntTerm}  
\end{eqnarray}

The length of $\gamma_{x_0,X}^\T$ restricted to the $j^\mathrm{th}$
subinterval equals $d(x_{j-1},x_j)$.
Since the $j^\mathrm{th}$ subinterval has length $t_j$ we have
$|\dot{\gamma}_{x_0,X}^\T|=d(x_{j-1},x_j)/t_j$ on the $j^\mathrm{th}$
subinterval. 
Hence the energy of $\gamma_{x_0,X}^\T$ on the $j^\mathrm{th}$ subinterval is
$\frac12\cdot t_j\cdot(d(x_{j-1},x_j)/t_j)^2 = \frac12\cdot
d(x_{j-1},x_j)^2/t_j$. 
Therefore
\begin{equation}
\E(\gamma_{x_0,X}^\T) = \frac12 \sum_{j=1}^r \frac{d(x_{j-1},x_j)^2}{t_j} .
\label{eq:EnTerm}
\end{equation}
Plugging (\ref{eq:IntTerm}), (\ref{eq:EnTerm}), and (\ref{eq:ZTerm})
 into (\ref{eq:Expand1}) we get
\begin{eqnarray}
\lefteqn{\hW_{t_{1}^n}\cdots \hW_{t_{r_n}^n}u(x_0)}\nonumber \\
&=&
\frac{1}{Z(\T_n,m)} \,
\int_{\M^{r_n}} 
\exp\left(-\frac12\E(\gamma_{x_0,X}^{\T_n})\right) 
\cdot\exp\left(\int_0^{\Laenge(\T_n)}\frac13 \scal(\gamma_{x_0,X}^{\T_n}(s))\, 
ds\right)\times\nonumber \\ 
&&
\times\tau(\gamma_{x_0,X}^{\T_n},\nabla)_{\Laenge(\T_n)}^{0}\times\nonumber \\ 
&&
\times
\prod_{j=1}^{r_n}
\exp\left(-\int_{\sigma_{j-1}(\T_n)}^{\sigma_j(\T_n)}
  \tau(\gamma_{x_0,X}^{\T_n},\nabla)_s^{\Laenge(\T_n)}\cdot
  V(\gamma_{x_0,X}^{\T_n}(s)) 
  \cdot\tau(\gamma_{x_0,X}^{\T_n},\nabla)_{\Laenge(\T_n)}^{s}\,ds\right)\cdot 
u(x_{r_n})\,dX .
\label{eq:Expand2}
\end{eqnarray}
Rewriting the integral over $\M^{r_n}$ as an integral over $\P_{x_0}(M,\T_n)$
finishes the proof of Theorem~\ref{thm:PathIntegral1}.  
\hfill$\Box$

\section{Alternative versions of the main theorem}
\label{sec:Alternative}

It is not mandatory to perform all six modifications in the proof of
Theorem~\ref{thm:PathIntegral1}.
If we cancel the sixth modification we obtain

\begin{thm}\label{thm:PathIntegral2}
Let $M$ be an $m$-dimensional closed Riemannian manifold, let $E$ be a vector
bundle over $M$ with a metric and a compatible connection $\nabla$.
Let $H = \nabla^*\nabla + V$ be a self-adjoint generalized Laplace
operator acting on sections in $E$.
Let $t>0$.

Then for any sequence of partitions $\T_n =(t_1^n, \ldots,t_{r_n}^n)$ with
$|\T_n| \to 0$ and $\Laenge(\T_n) \to t$ as $n \to \infty$ and for any $u\in
C^0(M,E)$ 
\begin{equation*}
\frac{1}{Z(\T_n,m)} \,
\int_{\P_x(M,\T_n)} 
\exp\left(-\frac12\E(\gamma)+\frac13\int_0^{\Laenge(\T_n)} \scal(\gamma(s))\, 
ds\right)\cdot \chi(\gamma,\T_n)\cdot
\tau(\gamma,\nabla)_{\Laenge(\T_n)}^{0}\times
\end{equation*}
\begin{equation*}
\times
\prod_{j=1}^{r_n}
\exp\left(-\int_{\sigma_{j-1}(\T_n)}^{\sigma_j(\T_n)}
  \tau(\gamma,\nabla)_s^{\Laenge(\T_n)}\cdot
  V(\gamma(s)) 
  \cdot\tau(\gamma,\nabla)_{\Laenge(\T_n)}^{s}\,ds\right)\cdot 
u(\gamma(\Laenge(\T_n)))\,\Dg
\end{equation*}
\begin{equation*}
\quad\quad\xrightarrow{n \to \infty} \quad\quad e^{-tH}u(x)
\end{equation*}
converges uniformly in $x$.
\hfill $\Box$
\end{thm}

Here the additional cut-off term is defined as 
$$
\chi(\gamma,\T) := 
\prod_{j=1}^r \chi(d(\gamma(\sigma_{j-1}(\T)),\gamma(\sigma_{j}(\T)))^2).
$$
Heuristically, this tells us that geodesic polygons having at least one ``long
edge'' do not contribute to the path integral.

Using Remark~\ref{rem:scal} we can generalize Theorem~\ref{thm:PathIntegral2}
to
\begin{thm}\label{thm:PathIntegral3}
Let $M$ be an $m$-dimensional closed Riemannian manifold, let $E$ be a vector
bundle over $M$ with a metric and a compatible connection $\nabla$.
Let $H = \nabla^*\nabla + V$ be a self-adjoint generalized Laplace
operator acting on sections in $E$.
We fix $\Lambda\in\R$.
Let $t>0$.

Then for any sequence of partitions $\T_n =(t_1^n, \ldots,t_{r_n}^n)$ with
$|\T_n| \to 0$ and $\Laenge(\T_n) \to t$ as $n \to \infty$ and for any $u\in
C^0(M,E)$ 
\begin{equation*}
\frac{1}{Z(\T_n,m)} \,
\int_{\P_x(M,\T_n)} 
\exp\left(-\frac12\E(\gamma)+\frac{\Lambda+1}{6}\int_0^{\Laenge(\T_n)}
\scal(\gamma(s))\,  
ds\right)\cdot \mu(\gamma,\T_n)^{\frac{\Lambda-1}{2}}\cdot\chi(\gamma,\T_n)
\times
\end{equation*}
\begin{equation*}
\times
\tau(\gamma,\nabla)_{\Laenge(\T_n)}^{0}
\prod_{j=1}^{r_n}
\exp\left(-\int_{\sigma_{j-1}(\T_n)}^{\sigma_j(\T_n)}
  \tau(\gamma,\nabla)_s^{\Laenge(\T_n)}
  V(\gamma(s)) 
  \tau(\gamma,\nabla)_{\Laenge(\T_n)}^{s}\,ds\right)
u(\gamma(\Laenge(\T_n)))\,\Dg
\end{equation*}
\begin{equation*}
\quad\quad\xrightarrow{n \to \infty} \quad\quad e^{-tH}u(x)
\end{equation*}
converges uniformly in $x$.
\hfill $\Box$
\end{thm}

Here, of course, the additional measure term is defined as 
$$
\mu(\gamma,\T) := 
\prod_{j=1}^r \mu(d(\gamma(\sigma_{j-1}(\T)),\gamma(\sigma_{j}(\T)))^2).
$$
The choice $\Lambda=-1$ is particularly interesting.
In this case the scalar curvature term disappears and $\mu(\gamma,\T)^{-1}\Dg$
is the product measure of {\em Euclidean} volume measures induced by
Riemannian normal coordinates.
Let us denote this measure by $\mathrm{D}\gamma$.
Then we get

\begin{cor}\label{cor:PathIntegral3}
Let $M$ be an $m$-dimensional closed Riemannian manifold, let $E$ be a vector
bundle over $M$ with a metric and a compatible connection $\nabla$.
Let $H = \nabla^*\nabla + V$ be a self-adjoint generalized Laplace
operator acting on sections in $E$.
Let $t>0$.

Then for any sequence of partitions $\T_n =(t_1^n, \ldots,t_{r_n}^n)$ with
$|\T_n| \to 0$ and $\Laenge(\T_n) \to t$ as $n \to \infty$ and for any $u\in
C^0(M,E)$ 
\begin{equation*}
\frac{1}{Z(\T_n,m)} \,
\int_{\P_x(M,\T_n)} 
\exp\left(-\frac12\E(\gamma)\right)\cdot\chi(\gamma,\T_n)
\cdot\tau(\gamma,\nabla)_{\Laenge(\T_n)}^{0}
\times
\end{equation*}
\begin{equation*}
\times
\prod_{j=1}^{r_n}
\exp\left(-\int_{\sigma_{j-1}(\T_n)}^{\sigma_j(\T_n)}
  \tau(\gamma,\nabla)_s^{\Laenge(\T_n)}
  V(\gamma(s)) 
  \tau(\gamma,\nabla)_{\Laenge(\T_n)}^{s}\,ds\right)
u(\gamma(\Laenge(\T_n)))\,\mathrm{D}\gamma
\end{equation*}
\begin{equation*}
\quad\quad\xrightarrow{n \to \infty} \quad\quad e^{-tH}u(x)
\end{equation*}
converges uniformly in $x$.
\hfill $\Box$
\end{cor}

The choices $\Lambda=1$ and $\Lambda=-1$ and hence $\Dg$ and $\mathrm{D}\gamma$
correspond to the measures on $\P_x(M,\T)$ induced by the $L^2$-metric and the
$H^1$-metric in \cite{AD} respectively.
Theorem~\ref{thm:PathIntegral3} gives an interpolation between these two cases.

Finally, since
$$
\left| \frac{t}{6}\scal(x) + \frac{t}{6}\scal(y) 
- \int_0^t \frac13 \scal(\gotxy(s))\,ds\right|
\quad\leq\quad 
C \cdot t \cdot d(x,y)
$$
Lemma~\ref{lem:KernKriterium} applies with $\alpha=\beta=1$ and shows that we
may replace $\int_0^t \frac13 \scal(\gotxy(s))\,ds$ by $\frac{t}{6}\scal(x)+
\frac{t}{6}\scal(y)$.
This shows

\begin{thm}\label{thm:PathIntegral4}
Let $M$ be an $m$-dimensional closed Riemannian manifold, let $E$ be a vector
bundle over $M$ with a metric and a compatible connection $\nabla$.
Let $H = \nabla^*\nabla + V$ be a self-adjoint generalized Laplace
operator acting on sections in $E$.
Let $t>0$.

Then for any sequence of partitions $\T_n =(t_1^n, \ldots,t_{r_n}^n)$ with
$|\T_n| \to 0$ and $\Laenge(\T_n) \to t$ as $n \to \infty$ and for any $u\in
C^0(M,E)$ 
\begin{equation*}
\frac{1}{Z(\T_n,m)} \,
\int_{\P_x(M,\T_n)} 
\exp\left(-\frac12\E(\gamma)\, ds\right)\cdot
\tau(\gamma,\nabla)_{\Laenge(\T_n)}^{0}\cdot
\prod_{j=1}^{r_n}
\exp\bigg(\frac{t_j}{6}\scal(\gamma(\sigma_{j-1}(\T_n))) +
\end{equation*}
\begin{equation*}
+\frac{t_j}{6}\scal(\gamma(\sigma_{j}(\T_n)))
-\int_{\sigma_{j-1}(\T_n)}^{\sigma_j(\T_n)}
  \tau(\gamma,\nabla)_s^{\Laenge(\T_n)}\cdot
  V(\gamma(s)) 
  \cdot\tau(\gamma,\nabla)_{\Laenge(\T_n)}^{s}\,ds\bigg)\cdot 
u(\gamma(\Laenge(\T_n)))\,\Dg
\end{equation*}
\begin{equation*}
\quad\quad\xrightarrow{n \to \infty} \quad\quad e^{-tH}u(x)
\end{equation*}
converges uniformly in $x$.
\end{thm}

This generalizes Corollary~{1.9} in \cite{AD}.
Note that the factor of $\frac16$ in front of the scalar curvature has to be
replaced by $\frac1{12}$ in \cite{AD} because there one considers the
semigroup $e^{-\frac{t}2\Delta}$.

\section{Approximation of the heat kernel}

Theorem~\ref{thm:PathIntegral1} is an approximation result for the heat
semigroup $e^{-tH}$ in the strong operator topology.
We will now refine this by uniformly approximating the integral kernel
$k_t(x,y)$ of $e^{-tH}$. 
According to the formula in Theorem~\ref{thm:PathIntegral1} we define for any
$x,y\in M$ and any partition $\T=(t_1,\ldots,t_r)$
\begin{equation*}
k_{\T}(x,y) :=\frac{1}{Z(\T,m)} \,
\int_{\P_x^y(M,\T)} 
\exp\left(-\frac12\E(\gamma)+\frac13\int_0^{\Laenge(\T)} \scal(\gamma(s))\, 
ds\right)\cdot
\tau(\gamma,\nabla)_{\Laenge(\T)}^{0}\times
\end{equation*}
\begin{equation}
\times
\prod_{j=1}^{r}
\exp\left(-\int_{\sigma_{j-1}(\T)}^{\sigma_j(\T)}
  \tau(\gamma,\nabla)_s^{\Laenge(\T)}\cdot
  V(\gamma(s)) 
  \cdot\tau(\gamma,\nabla)_{\Laenge(\T)}^{s}\,ds\right)\,\Dg .
\label{eq:DefApproxKern}
\end{equation}
Now the statement of Theorem~\ref{thm:PathIntegral1} is
$$
\int_M k_{\T_n}(x,y)u(y)dy 
\xrightarrow{n\to \infty}
\int_M k_t(x,y)u(y)dy 
$$
uniformly in $x$ provided $|\T_n|\to 0$, $\Laenge(\T_n)\to t$, and $u\in C^0(M,E)$.
We will improve this to the statement that $k_{\T_n} \to k_t$ in
$C^0(M\times M,E\boxtimes E^*)$, at least for suitable sequences of partitions
$\T_n$.

\begin{thm}\label{thm:PathKernel1}
Let $M$ be an $m$-dimensional closed Riemannian manifold, let $E$ be a vector
bundle over $M$ with a metric and a compatible connection $\nabla$.
Let $H = \nabla^*\nabla + V$ be a self-adjoint generalized Laplace
operator acting on sections in $E$.
Denote the heat kernel of $H$ by $k_t(x,y)$.

Given $t>0$ there exists a sequence of partitions $\T_n$ with $\Laenge(\T_n)=t$ and
$|\T_n| \to 0$ such that 
$$
\| k_{\T_n} - k_t\|_{C^0(M \times M)} \xrightarrow{n\to\infty} 0
$$
where $k_{\T_n}$ are defined in (\ref{eq:DefApproxKern}).
\end{thm}

\begin{proof}
We fix $t>0$ and $\eps>0$.
Given $y\in M$ the section
\begin{eqnarray*}
x \mapsto k_{(t_r)}(x,y) 
&=& 
(4\pi t_r)^{-m/2}\exp\left(-\frac{d(x,y)^2}{4t_r}\right)
\cdot \tau(\gamma,\nabla)_{t}^{t-t_r} \times\\
&&
\times\exp\int_{t-t_r}^t \left(\frac13\scal(\gamma) -
\tau(\gamma,\nabla)_s^{t}\cdot V(\gamma(s)) 
  \cdot\tau(\gamma,\nabla)_{t}^{s}\right)ds,
\end{eqnarray*}
$\gamma=\gamma_{x,y}^{t-t_r,t}$, is a Gaussian approximation to the delta
function at $y$ for $t_r$ small.
Hence $k_{(t_r)}(\cdot,y) \to \delta_y$ in the sense of
distributions as $t_r \searrow 0$.
Note that $\delta_y$ is a distributional section in the bundle $E \otimes
E_y^*$.
Also note that the term $\exp\int_{t-t_r}^t \left(\frac13\scal(\gamma) -
\tau(\gamma,\nabla)_s^{t}\cdot V(\gamma(s)) \cdot\tau(\gamma,\nabla)_{t}^{s}
\right)ds$ is irrelevant here because its difference with the identity is of
order $\OOOO(t_r)$ in the operator norm.
Since $e^{-sH}$ is smoothing we have
$$
\|k_s(\cdot,y) - e^{-sH}k_{(t_r)}(\cdot,y) \|_{C^0(M)}
\quad=\quad
\|e^{-sH}(\delta_y - k_{(t_r)}(\cdot,y)) \|_{C^0(M)}
\quad\leq\quad \eps
$$
if $t_r$ is sufficiently small and $s$ is contained in any fixed compact
subinterval of $(0,\infty)$.
In particular, we get with $s=t-t_r$
\begin{equation}\label{eq:KernAp1}
\| k_{t-t_r}(\cdot,y) - e^{-(t-t_r)H}k_{(t_r)}(\cdot,y) \|_{C^0(M)} 
\quad\leq\quad \eps.
\end{equation}
Since the heat kernel of $H$ is continuous we have for sufficiently small
$t_r$
\begin{equation}\label{eq:KernAp2}
\| k_t(\cdot,y) - k_{t-t_r}(\cdot,y) \|_{C^0(M)} \quad\leq\quad \eps.
\end{equation}
By Theorem~\ref{thm:PathIntegral1} with $u=k_{(t_r)}(\cdot,y)$ we have for
sufficiently fine partitions $\T'=(t_1,\ldots,t_{r-1})$ of $[0,t-t_r]$ that 
\begin{eqnarray}
\eps 
&\geq&
\left\|e^{-(t-t_r)H}k_{(t_r)}(\cdot,y) - \int_M k_{\T'}(\cdot,z)k_{(t_r)}
(z,y) \, dz \right\|_{C^0(M)}\nonumber \\ 
&=&
\left\|e^{-(t-t_r)H}k_{(t_r)}(\cdot,y) - k_{\T}(\cdot,y)\right\|_{C^0(M)}
\label{eq:KernAp3}
\end{eqnarray}
with $\T=(t_1,\ldots,t_{r-1},t_r)$.
Combining (\ref{eq:KernAp1}), (\ref{eq:KernAp2}), and (\ref{eq:KernAp3}) we get
$$
\| k_t(\cdot,y) - k_{\T}(\cdot,y) \|_{C^0(M)} \quad\leq\quad 3\eps.
$$
By compactness of $M$ the estimates can be made uniformly in $y$ (compare
Lemma~\ref{lem:compactuniform}) so that
$$
\| k_t - k_{\T} \|_{C^0(M\times M)} \quad\leq\quad 3\eps.
$$
\end{proof}

\begin{rem}\label{rem:PathKernel1}
The proof shows that the sequence of partitions $\T_n$ can be found as
follows.
First choose $t_r$ small enough ($r$ itself is not specified yet), then choose
$\T'=(t_1,\ldots,t_{r-1})$ a sufficiently fine partition of $[0,t-t_r]$ (now
$r$ is determined), then $\T=(\T',t_r)$ yields a good approximation for the
heat kernel.

Therefore it is possible to choose the same sequence of partitions $\T_n$
working simultaneously for finitely many given operators $H_1,\ldots,H_k$
(possibly acting on sections in different bundles).
\end{rem}

Uniform convergence of integral kernels implies convergence of the integral
operators as operators from $L^2$ to $C^0$.
Hence we have

\begin{cor}\label{cor:PathKernel}
Let $M$ be an $m$-dimensional closed Riemannian manifold, let $E$ be a vector
bundle over $M$ with a metric and a compatible connection $\nabla$.
Let $H = \nabla^*\nabla + V$ be a self-adjoint generalized Laplace
operator acting on sections in $E$.
Let $t>0$.

Then there exists sequence of partitions $\T_n =(t_1^n, \ldots,t_{r_n}^n)$ with
$|\T_n| \to 0$ and $\Laenge(\T_n) = t$ such that for any $u\in L^2(M,E)$ 
\begin{equation*}
\frac{1}{Z(\T_n,m)} \,
\int_{\P_x(M,\T_n)} 
\exp\left(-\frac12\E(\gamma)+\frac13\int_0^{t} \scal(\gamma(s))\, 
ds\right)\cdot
\tau(\gamma,\nabla)_{t}^{0}\times
\end{equation*}
\begin{equation*}
\times
\prod_{j=1}^{r_n}
\exp\left(-\int_{\sigma_{j-1}(\T_n)}^{\sigma_j(\T_n)}
  \tau(\gamma,\nabla)_s^{t}\cdot
  V(\gamma(s)) 
  \cdot\tau(\gamma,\nabla)_{t}^{s}\,ds\right)\cdot 
u(\gamma(t))\,\Dg
\end{equation*}
\begin{equation*}
\quad\quad\xrightarrow{n \to \infty} \quad\quad e^{-tH}u(x)
\end{equation*}
converges uniformly in $x$.
\hfill$\Box$
\end{cor}

\begin{rem}
In contrast to Theorem~\ref{thm:PathIntegral1} where {\em all} sequences of
partitions $\T_n$ with $\Laenge(\T_n) \to t$ and $|\T_n|\to 0$ yield the right
approximation Theorem~\ref{thm:PathKernel1} and
Corollary~\ref{cor:PathKernel} make statements only for {\em some} such
sequences.
We do not know whether Theorem~\ref{thm:PathKernel1} and
Corollary~\ref{cor:PathKernel} still hold for all such sequences.
For the application in the next section however this difference is irrelevant.
\end{rem}

\begin{rem}
One may use one of the alternative versions of the main theorem in
Section~\ref{sec:Alternative} instead of Theorem~\ref{thm:PathIntegral1}.
For example, using Theorem~\ref{thm:PathIntegral2} we get
Theorem~\ref{thm:PathKernel1} where the definition in (\ref{eq:DefApproxKern})
is replaced by
\begin{equation*}
k_{\T}(x,y) :=\frac{1}{Z(\T,m)} \,
\int_{\P_x^y(M,\T)} 
\exp\left(-\frac12\E(\gamma)+\frac13\int_0^{\Laenge(\T)} \scal(\gamma(s))\, 
ds\right)\cdot
\tau(\gamma,\nabla)_{\Laenge(\T)}^{0}\times
\end{equation*}
\begin{equation}
\quad\times
\chi(\gamma,\T)\cdot\prod_{j=1}^{r}
\exp\left(-\int_{\sigma_{j-1}(\T)}^{\sigma_j(\T)}
  \tau(\gamma,\nabla)_s^{\Laenge(\T)}\cdot
  V(\gamma(s)) 
  \cdot\tau(\gamma,\nabla)_{\Laenge(\T)}^{s}\,ds\right)\,\Dg .
\end{equation}
\end{rem}

\section{The Hess-Schrader-Uhlenbrock estimate}

Let $H=\nabla^*\nabla +V$ be a self-adjoint Laplace type operator acting on
sections in a Riemannian or Hermitian vector bundle $E$ over a closed
Riemannian manifold $M$. 
Let $v\in C^\infty(M)$ be a real valued function such that $V\geq v$
everywhere, i.~e., for each $x\in M$ the eigenvalues of $V(x)$ are bounded
below by $v(x)$. 
Let $\tilde H$ be the ``comparison operator'' $\tilde H = \Delta + v$ acting
on functions.
Let $k_t$ and $\tilde k_t$ be the corresponding heat kernels.

The main result of \cite{HSU} is the estimate
\begin{equation} \label{eq:Uhlenbrock}
|k_t(x,y)| \leq \tilde k_t(x,y)
\end{equation}
for all $x,y\in M$ and $t>0$.
This directly implies 
$$
\Tr(e^{-tH}) \leq \Tr(e^{-t\tilde H})
$$
for all $t>0$ and hence $\lambda_1(H) \geq \lambda_1(\tilde H)$ where
$\lambda_1$ denotes the smallest eigenvalue.
For further applications see e.~g.\ \cite{B,Gr}.

The proof in \cite{HSU} is based on a Kato inequality.
We demonstrate here that (\ref{eq:Uhlenbrock}) is a direct consequence of our
path integral approximation for the heat kernel.
Namely, for any partition $\T$ we have
\begin{eqnarray*}
|k_{\T}(x,y)|
&\leq&
\frac{1}{Z(\T,m)} \,
\int_{\P_x^y(M,\T)} 
\exp\left(-\frac12\E(\gamma)+\frac13\int_0^{\Laenge(\T)} \scal(\gamma(s))\, 
ds\right)\cdot
|\tau(\gamma,\nabla)_{\Laenge(\T)}^{0}|\times\\
&&
\times
\prod_{j=1}^{r}
\left|\exp\left(-\int_{\sigma_{j-1}(\T)}^{\sigma_j(\T)}
  \tau(\gamma,\nabla)_s^{\Laenge(\T)}\cdot
  V(\gamma(s)) 
  \cdot\tau(\gamma,\nabla)_{\Laenge(\T)}^{s}\,ds\right)\right|\,\Dg\\
&\leq&
\frac{1}{Z(\T,m)} \,
\int_{\P_x^y(M,\T)} 
\exp\left(-\frac12\E(\gamma)+\frac13\int_0^{\Laenge(\T)} \scal(\gamma(s))\, 
ds\right)\cdot
1\times\\
&&
\times
\prod_{j=1}^{r}
\exp\left(-\int_{\sigma_{j-1}(\T)}^{\sigma_j(\T)}  v(\gamma(s)) 
  \,ds\right)\,\Dg\\
&=&
\tilde k_\T(x,y).
\end{eqnarray*}
Now (\ref{eq:Uhlenbrock}) follows directly from Theorem~\ref{thm:PathKernel1}
and Remark~\ref{rem:PathKernel1}.

\section{A trace formula}

Parallel transport (with respect to a connection $\nabla$) along a closed loop
$\gamma$ is called its {\em holonomy}.
Note that the holonomy $\hol(\gamma,\nabla)$ is depends on the base point of
the loop only up to conjugation.
Therefore $\tr(\hol(\gamma,\nabla))$ is independent of the base point.

\begin{thm}\label{thm:Trace}
Let $M$ be an $m$-dimensional closed Riemannian manifold, let $E$ be a vector
bundle over $M$ with a metric and a compatible connection $\nabla$.
Let $H = \nabla^*\nabla + V$ be a self-adjoint generalized Laplace
operator acting on sections in $E$.
Let the potential $V$ be scalar.

Given $t>0$ there exists a sequence of partitions $\T_n$ with $\Laenge(\T_n)=t$ and
$|\T_n| \to 0$ such that 
\begin{eqnarray*}
\lefteqn{\Tr(e^{-tH})}\\
&=&
\lim_{n\to\infty}
\frac{1}{Z(\T_n,m)} \,
\int_{\Pcl(M,\T_n)} 
\exp\left(-\frac12\E(\gamma)
+\int_0^{t}\left(\frac13 \scal(\gamma(s))-V(\gamma(s))\right)\, 
ds\right)\times\\
&&
\quad\quad\times
\tr(\hol(\gamma,\nabla))\,\Dg .
\end{eqnarray*}
\end{thm}

\begin{proof}
Let $\T_n$ be a sequence of partitions as in Theorem~\ref{thm:PathKernel1}.
Then the pointwise trace $\tr(k_{\T_n})$ converges uniformly to
$\tr(k_t)$  and we have
\begin{eqnarray*}
\Tr(e^{-tH})
&=&
\int_M \tr(k_t(x,x))\, dx\\
&=&
\lim_{n\to\infty} \int_M \tr(k_{\T_n}(x,x))\, dx \\
&=&
\lim_{n\to\infty} \int_M \frac{1}{Z(\T_n,m)} \int_{\P_x^x(M,\T_n)}
\exp(\cdots)\cdot \tr(\tau(\gamma,\nabla))\,\Dg\, dx\\
&=&
\lim_{n\to\infty} \frac{1}{Z(\T_n,m)} \int_{\Pcl(M,\T_n)}
\exp(\cdots)\cdot \tr(\hol(\gamma,\nabla))\,\Dg .
\end{eqnarray*}
\end{proof}

Of course, one also obtains a trace formula if $V$ is not scalar.
In this case $V$ and $\tau$ cannot be separated and instead of 
$\exp\left(-\int_0^tV(\gamma(s))ds\right)\tr(\hol(\gamma,\nabla))$ one gets
$\tr\left(\tau(\gamma,\nabla)_{t}^{0}\cdot
\prod_{j=1}^{r_n}
\exp\left(-\int_{\sigma_{j-1}(\T_n)}^{\sigma_j(\T_n)}
  \tau(\gamma,\nabla)_s^{t}\cdot
  V(\gamma(s)) 
  \cdot\tau(\gamma,\nabla)_{t}^{s}\,ds\right)\right)$.

\appendix

\section{A Gaussian estimate}

For the sake of completeness we include the following technical result.

\begin{lem}\label{lem:GaussEst}
Let $B: \R^m \times \R^m \to \R$ be a symmetric bilinear form.
Let $f:[0,t_0] \times \R^m \to \C^k$ be a $C^1$-map with finite
$C^1(\R^m)$-norm, i.~e., there exists a constant $C>0$ such that $|f(t,\xi)|
\leq C$ and $|\abl{\xi_j}f(t,\xi)|\leq C$ for all $(t,\xi) \in [0,t_0] \times
\R^m$. 

Then, as $t \searrow 0$,
$$
\int_{\R^m}\frac{e^{-\frac{|\xi|^2}{4t}}}{(4\pi t)^{m/2}} \cdot B(\xi,\xi)
\cdot f(t,\xi)\, d\xi
=
2t\cdot \tr(B) \cdot \int_{\R^m}\frac{e^{-\frac{|\xi|^2}{4t}}}{(4\pi t)^{m/2}}
\cdot f(t,\xi)\, d\xi + \OOOO(t^{3/2}) .
$$
The constant in the $\OOOO(t^{3/2})$-term depends only on $C$, $m$, and an upper
bound on $|B|$.
\end{lem}

\begin{proof}
Without loss of generality we assume that the Cartesian coordinates are chosen
such that $B$ is diagonalized, $B(\xi,\eta) = \sum_{j=1}^m \lambda_j
\xi_j\eta_j$.
Fix $t>0$.
The smooth vector field 
$$
X:= \frac{e^{-\frac{|\xi|^2}{4t}}}{(4\pi t)^{m/2}} \cdot 
\sum_{j=1}^m \lambda_j \xi_j \abl{\xi_j}
$$
is rapidly decreasing as $|\xi| \to \infty$.
Hence we may integrate by parts to get
\begin{eqnarray*}
\int_{\R^m} \div\left(X\right) \cdot f(t,\xi)\, d\xi
&=&
- \int_{\R^m} \< X, \grad  f(t,\xi)\> \, d\xi\\
&=&  
- \int_{\R^m} \frac{e^{-\frac{|\xi|^2}{4t}}}{(4\pi t)^{m/2}} \cdot 
\sum_{j=1}^m \lambda_j \xi_j \abl{\xi_j} f(t,\xi) \, d\xi .
\end{eqnarray*}
Thus 
\begin{equation*}
\left| \int_{\R^m} \div\left(X\right) \cdot f(t,\xi)\, d\xi \right|
\quad\leq\quad
C_1 \cdot \int_{\R^m} \frac{e^{-\frac{|\xi|^2}{4t}}}{(4\pi t)^{m/2}} \cdot
|\xi|\, d\xi .
\end{equation*}
Using $\tau < \exp(\tau^2)$ for all $\tau\in\R$ with $\tau =
\frac{|\xi|}{\sqrt{8t}}$ yields $|\xi| \leq e^{\frac{|\xi|^2}{8t}} \cdot
\sqrt{8t}$ and hence
\begin{equation}
\label{eq:DivAbsch}
\left| \int_{\R^m} \div\left(X\right) \cdot f(t,\xi)\, d\xi \right|
\quad\leq\quad
C_1 \cdot \sqrt{8t}\cdot \int_{\R^m} \frac{e^{-\frac{|\xi|^2}{8t}}}{(4\pi
  t)^{m/2}} \, d\xi 
\quad=\quad
C_2 \cdot \sqrt{t} .
\end{equation}
On the other hand,
\begin{eqnarray*}
\div(X)
&=&
\sum_{j=1}^m \abl{\xi_j}\left(\frac{e^{-\frac{|\xi|^2}{4t}}}{(4\pi t)^{m/2}} \cdot 
\lambda_j \xi_j\right)\\
&=&
\sum_{j=1}^m \lambda_j \left(-\frac{\xi_j^2}{2t} +1\right)
\frac{e^{-\frac{|\xi|^2}{4t}}}{(4\pi t)^{m/2}}\\
&=&
\left(-\frac{B(\xi,\xi)}{2t} +\tr(B)\right)
\frac{e^{-\frac{|\xi|^2}{4t}}}{(4\pi t)^{m/2}} .
\end{eqnarray*}
Multiplication of this equation with $2t f(t,\xi)$ and integration over $\R^m$
together with (\ref{eq:DivAbsch}) yields the assertion.
\end{proof}

\end{document}